\documentclass[11pt]{article}
\usepackage{enumerate}
\usepackage{mathtools}
\usepackage{amssymb}
\usepackage{amsmath}
\usepackage{amsthm}
\usepackage{mathrsfs}
\usepackage{graphicx}
\usepackage{multirow}

\usepackage{float}
\usepackage{subfig}

\usepackage{hyperref}

\usepackage{cleveref}

\newtheorem{theorem}{Theorem}

\newtheorem{example}{Example}

\usepackage{bm}
\usepackage{tabularx}
\usepackage{color}
\usepackage{multirow}
\usepackage{algorithm}
\usepackage{algorithmic}
\usepackage{threeparttable}
\usepackage{listings}

\def\bb{\pmb{b}}
\def\be{\pmb{e}}
\def\bx{\pmb{x}}
\def\by{\pmb{y}}

\allowdisplaybreaks

\title{Compensated Splitting For Generalized Lyapunov Equations}
\author{Hongjia Chen%
\thanks{Department of Mathematics, Nanchang University, 999 Xuefu Road, Nanchang, Jiangxi, 330031, China.
        Supported by the National Natural Science Foundation of China under grant No.~12371379 and No.~12361080,
and Jiangxi Provincial Natural Science Foundation under grant No.~20224BAB211006.
        Email: {\tt chenhongjia@ncu.edu.cn}.}
\and Ren-Cang Li%
\thanks{Department of Mathematics, University of Texas at Arlington, Arlington, TX 76019-0408, USA.
        Supported in part by NSF DMS-2407692.
        Email: {\tt rcli@uta.edu}.}
}

\DeclareMathOperator{\Id}{Id}

\DeclareMathOperator{\vex}{vec}

\def\bbC{\mathbb{C}}
\def\bbR{\mathbb{R}}

\def\scrA{\mathscr{A}}
\def\scrM{\mathscr{M}}
\def\scrN{\mathscr{N}}

\def\wtd{\widetilde}
\def\bE{\bar{E}}
\def\bN{\bar{N}}

\DeclareMathOperator{\F}{F}
\DeclareMathOperator{\HH}{H}
\DeclareMathOperator{\T}{T}
\DeclareMathOperator{\tr}{tr}
\DeclareMathOperator{\NRes}{NRes}

\numberwithin{algorithm}{section}
\numberwithin{equation}{section}
\numberwithin{figure}{section}
\numberwithin{table}{section}
\numberwithin{theorem}{section}

\begin{document}
\maketitle

\begin{abstract}
The generalized Lyapunov equation has a natural splitting that leads to the standard fixed-point iteration
(sFPI) for its numerical solution. sFPI is convergent, for an arbitrarily given initial guess,
if and only if the spectral radius of the associated linear operator
is smaller than $1$. This means that sFPI may diverge  if the spectral radius is $1$ or bigger.
In this paper, we propose a compensated splitting scheme
that can hopefully reduce the spectral radius
so that the resulting compensated fixed-point iteration (cFPI) has a better convergence property, namely,
cFPI may still converge even if sFPI does not or cFPI converges faster than
sFPI does when the latter is also convergent.
Numerical results are presented to demonstrate the superiority of cFPI to sFPI.


\end{abstract}

\section{Introduction}\label{sec:intro}
Our focus in this paper is the numerical solution to the generalized Lyapunov equation
\begin{equation}\label{eq:genlyap}
AX + XA^{\HH} + \sum_{j=1}^{m}N_j X N_j^{\HH} + BB^{\HH} = 0,
\end{equation}
where the matrices $A, N_j\in\bbC^{n\times n}$ and $B\in\bbC^{n\times \ell}$ (usually $\ell\ll n$).
Such an equation often arises
from the model order reduction of bilinear control systems and linear parameter-varying systems,
the stability analysis of linear stochastic differential
equations, as well as the linearization of stochastic continuous-time
algebraic Riccati equations~\cite{benner2013low,beda:2011,damm2008direct,damm2001newton,hukl:2025}.
The constant term  in \eqref{eq:genlyap} takes such a special form $BB^{\HH}$  because of
the sources of the generalized Lyapunov equation in applications.
As far as the development in this paper is concerned, the constant term can be replaced with any square matrix $C$.
Imposing that $C$ is Hermitian will generally imply that the solution $X$ is Hermitian, too.
In most existing literature, all matrices in \eqref{eq:genlyap} are assumed to be real, but we let them be complex for generality
and that all matrices are real is simply a special case.


Theoretically,
\eqref{eq:genlyap} is an $n^2$-by-$n^2$ linear system of  equations and can be turned into the familiar standard form of a linear
system of  equations:
\begin{equation}\label{eq:GLyap2ST}
\scrA\vex(X)\equiv\Big[I_n\otimes A+\bar{A}\otimes I_n+\sum_{j=1}^m\bN_j\otimes N_j\Big]\vex(X)=-\vex(BB^{\HH}),
\end{equation}
where $\otimes$ is the matrix Kronecker product, $\bar{A}$ and $\bN_j$ denote the entrywise complex conjugate of $A$ and $N_j$, and $\vex(\cdot)$ vectorizes a matrix into a column vector by stacking up the columns of the matrix. It follows from this standard form that if
$\scrA\in\bbC^{n^2\times n^2}$ is nonsingular then \eqref{eq:genlyap} has a unique solution
and the solution $X$ is Hermitian, but it is not clear if $X$ is definite or not.
Another implication of \eqref{eq:GLyap2ST} is the temptation to solve the original
\eqref{eq:genlyap} via \eqref{eq:GLyap2ST}. Unfortunately, any direct solver such as variants of Gaussian elimination
and is not feasible unless $n$ is very small (a few tens) because of its $O(n^4)$ storage and $O(n^6)$ computational
complexity \cite{demm:1997,govl:2013}.
Usual iterative methods for linear systems \cite{demm:1997,gree:1997,saad:2003} are straightforwardly applicable without
the need to form $\scrA$ explicitly.

The generalized Lyapunov equation \eqref{eq:genlyap} has received increasing attention over the last decades.
The simplest numerical method is through the fixed-point technique, and it is based on  splitting the linear system \eqref{eq:genlyap} naturally into $\scrM(X)=\scrN(X)-BB^{\HH}$, where
\begin{equation}\label{eq:NatSplit}
\scrM(X):=AX + XA^{\HH}, \quad
\scrN(X)=-\sum_{j=1}^{m}N_j X N_j^{\HH}.
\end{equation}
Both $\scrM(\cdot)$ and $\scrN(\cdot)$ are linear operators
from $\bbC^{n\times n}$ to itself, and in particular, $\scrN(X)\preceq 0$ (negative semidefinite) for any $X\succeq 0$
(positive semidefinite). When $\scrM$ is invertible, such a splitting yields the following standard fixed-point
iteration (sFPI): given $X_0$, for $k=1,2,\ldots$ solve
\begin{equation}\label{eq:NatItn}
\scrM(X_k)\equiv AX_k+X_kA^{\HH}=\scrN(X_{k-1})-BB^{\HH}
\end{equation}
for $X_k$ until convergence. Suppose that $\scrM$ is invertible.
A sufficient condition for the convergence of sFPI  \eqref{eq:NatItn}  is
\begin{equation}\label{eq:nat-rho<1}
\rho(\scrM^{-1}\scrN)<1,
\end{equation}
where $\rho(\scrM^{-1}\scrN)$ is the spectral radius of $\scrM^{-1}\scrN$ as a linear operator in
$\bbC^{n\times n}$. The condition implies that necessarily
$$
\scrM-\scrN=\scrM\,(\Id-\scrM^{-1}\scrN)
$$
is invertible, where $\Id$ denotes the identity operator in $\bbC^{n\times n}$.
But,
if $\rho(\scrM^{-1}\scrN)\ge 1$, then the scheme is in general divergent.
Our goal in this paper is to rectify $\rho(\scrM^{-1}\scrN)\ge 1$ to yield a new iterative scheme that
may still converge. Roughly speaking, $\rho(\scrM^{-1}\scrN)\ge 1$ is caused by the $\scrN$-part being relative too ``heavy''.
Hence in order to overcome the violation, we have to take out a portion of the $\scrN$-part to compensate
the $\scrM$-part. For that purpose, we propose to somehow decompose the $\scrN$-part as
\begin{equation}\label{eq:scrN-split}
-\scrN(X)=EX+XE^{\HH}-\wtd\scrN(X)
\end{equation}
in such a way that  $\wtd\scrN(\cdot)$ is made relatively less ``heavy'' than before and then compensate
$\scrM(X)$ with $EX+XE^{\HH}$, leading to
\begin{equation}\label{eq:t-scrMscrN}
\scrM(X)-\scrN(X)=\underbrace{[\scrM(X)+EX+XE^{\HH}]}_{=:\wtd\scrM(X)}-\wtd\scrN(X),
\end{equation}
which we will call a  {\em compensated splitting}.
As a result, we can establish the following compensated fixed-point
iteration (cFPI): given $X_0$, for $k=1,2,\ldots$ solve
\begin{equation}\label{eq:NatItn-pert}
\wtd\scrM(X_k)\equiv (A+E)X_k+X_k(A+E)^{\HH}=\wtd\scrN(X_{k-1})-BB^{\HH}
\end{equation}
for $X_k$ until convergence.
An ideal compensated splitting \eqref{eq:t-scrMscrN} should minimize $\rho(\wtd\scrM^{-1}\wtd\scrN\,)$ subject to
$\wtd\scrM$ being invertible. Conceivably, seeking such an ideal compensated splitting is utterly impractical.
Instead, we seek somehow to ``minimize''
$\wtd\scrN(\cdot)$ optimally in certain sense.
In the end, we come up with an  ``optimal'' $E$ that takes a
very elegant form
\begin{equation}\label{eq:optimal(E):intro}
E=\frac{2}{2n+\mu}\sum_{j=1}^m\tr(\bN_{j})\left[N_j-\frac{2}{4n+\mu}\tr(N_{j})\cdot I_n\right],
\end{equation}
where $\mu\ge 0$ is a regularizing parameter to control
the magnitude of $E$.
How it is derived will be detailed in \cref{sec:splitting}, where precisely in what criterion this $E$
is optimal will also be made concrete.

FPI-type methods are simple to design and implement.
Beyond such methods as we have discussed so far, not surprisingly, modern iterative techniques
for linear systems (see, e.g., \cite{demm:1997,govl:2013,gree:1997,saad:2003})
have seen their extensions to the generalized Lyapunov equation \eqref{eq:genlyap}, including,
for example, the bilinear ADI method \cite{benner2013low},
the greedy low-rank methods \cite{krsi:2015},
Krylov subspace-type methods \cite{jamp:2018},
and exploitations of the fixed-point iteration \cite{shank2016efficient},
a residual-based generalized rational-Krylov-type subspace method \cite{breiten2019residual},
and a method based on the Sherman--Morrison--Woodbury formula
\cite{damm2008direct}.
However, in this paper we will be focusing solely on rescuing sFPI from its divergence or making it converge faster.
Any potential combination of the key idea of this paper with the state-of-the-art subspace-based projection techniques
will be investigated in forthcoming studies.


The rest of this paper is organized as follows.
In  \cref{sec:splitting}, a compensated splitting
is proposed and its main result is \eqref{eq:optimal(E):intro}.
Numerical experiments are reported in \cref{sec:egs}, demonstrating
the new method's superiority to sFPI. Finally, we draw conclusions and make comments in \cref{sec:conclude}.

{\bf Notation}.
$\bbR^{m\times n}$ and $\bbC^{m\times n}$ denote the sets of real
and complex $m\times n$ matrices, respectively;
$\bbR^n=\bbR^{n\times1}$ and $\bbC^n=\bbC^{n\times1}$, and $I_n\in\bbR^{n\times n}$ is the identity matrix.
For a matrix $X$, $X^{\T}$, $X^{\HH}$ and $\bar X$ denote its transpose, complex conjugate transpose
and entrywise complex conjugate, respectively.  $\|X\|_2$ and $\|X\|_{\F}$
are the spectral and Frobenius norms of a matrix $X$, respectively, and, if $X$ is square, $\tr(X)$ takes its trace.
$\Re(x)$ and $\Im(x)$ extract the real and imaginary part of a complex number $x$.
Other notation will be defined at its first appearance.

\section{Compensated Splitting }\label{sec:splitting}
As we discussed in section~\ref{sec:intro}, an ideal compensated splitting \eqref{eq:t-scrMscrN}
should make $\rho(\wtd\scrM^{-1}\wtd\scrN\,)$ as small as possible. Unfortunately, that is not numerically feasible, and
so we will have to be less ambitious.
We start by making $\wtd\scrN(X)$ as insignificant as possible over all $X\in\bbC^{n\times n}$. To that end,
we attempt to seek $E\in\bbC^{n\times n}$ by
\begin{equation}\label{eq:SARE-CS}
\min_E\Bigl\{\max_{\|\vex(X)\|_2\le 1}\Bigl\|\sum_{j=1}^m N_jXN_j^{\HH}-(EX+XE^{\HH})\Bigr\|_{\F}^2+\mu\|E\|_{\F}^2\Bigr\},
\end{equation}
where $\mu\|E\|_{\F}^2$ is a regularization term that directly controls the magnitude of $E$, and $\mu\ge 0$
is the regularizing parameter that has to be chosen properly.
Notice that
\begin{multline}\label{eq:SARE-CS'}
\max_{\|\vex(X)\|_2\le 1}\Bigl\|\Bigl(I_n\otimes E+\bar E\otimes I_n-\sum_{j=1}^{m}\bar N_j\otimes N_j\Bigr)\vex(X)\Bigr\|_2^2
        +\mu\|E\|_{\F}^2 \\
    =\Bigl\|I_n\otimes E+\bar E\otimes I_n-\sum_{j=1}^{m}\bar N_j\otimes N_j\Bigr\|_2^2
        +\mu\|E\|_{\F}^2.
\end{multline}
which involves the spectral norms of an $n^2\times n^2$ matrix, and conceivably
minimizing the expression in \eqref{eq:SARE-CS'} over $E$ is a difficult thing to do, if  at all possible.
Instead, we lower our ambition and solve a relaxed problem of
\eqref{eq:SARE-CS'}, namely,
\begin{equation}\label{eq:SARE-CS''}
\min_E\Bigl\{f(E):=\Bigl\|I_n\otimes E+\bar E\otimes I_n-\sum_{j=1}^{m}\bar N_j\otimes N_j\Bigr\|_{\F}^2
            +\mu\|E\|_{\F}^2\Bigr\}.
\end{equation}
It can be seen that the objective $f(E)$ in \eqref{eq:obj:real}
is a non-negative quadratic function in the real and imaginary parts of the entries of $E$ and it is convex.
Hence any stationary point of \eqref{eq:SARE-CS''} is a global minimum point. Our main result in this section is stated in the following theorem, where $\tr(\cdot)$ takes the trace of a square matrix.

\begin{theorem}\label{thm:GCS(complex)}
An optimal $E$ to \eqref{eq:SARE-CS''} is given by
\begin{equation}\label{eq:optimal(E):cpx}
E=\frac{2}{2n+\mu}\sum_{j=1}^m\tr(\bN_{j})\left[N_j-\frac{2}{4n+\mu}\tr(N_{j})\cdot I_n\right].
\end{equation}
Moreover, the minimizer is unique for the real case,  and for the complex case with $\mu>0$, but it is not unique for the complex case with $\mu=0$ in which case
this optimal $E$ is the one with the smallest
Frobenius norm.
\end{theorem}

We will defer the proof of this theorem to subsections~\ref{ssec:E-real-proof} and \ref{ssec:E-complex-proof} below. A couple of
immediately questions are:
1) what is the invertibility of the operator $\wtd\scrM$?
2) does $E$ in \eqref{eq:optimal(E):cpx} always decrease the spectral radius in question? The answers to both questions are negative,
as our later numerical experiments in \cref{sec:egs} suggest. The following identity may provide some understanding as to
when  $E$ in \eqref{eq:optimal(E):cpx} will work better. Notice that
\begin{equation}\label{eq:trE}
\tr(E)=\frac 2{4n+\mu}\sum_{j=1}^m|\tr(N_j)|^2\ge 0.
\end{equation}
Suppose that the magnitudes of all $N_j$ are small and then so is that of
$E$ in \eqref{eq:optimal(E):cpx}, and thus $E$ can be considered a perturbation.
Since
$$
\Re\tr(A+E)=\Re\tr(A)+\tr(E)\ge\Re\tr(A),
$$
$E$ will move some of the eigenvalues rightward.
Now if the eigenvalues of $A$ lie in the left-half plane, then collectively these eigenvalues will be moved rightward towards $0$,
making $\wtd\scrM$ ``less invertible'' than $\scrM$. In such a case there is a possibility that any gain in ``minimizing'' $\wtd\scrN$ may get offset by the
worsening in $\wtd\scrM$ along the way. On the other hand, if  the eigenvalues of $A$ lie in the right-half plane, then collectively
the eigenvalues of $A$ will also be moved rightward but further away from $0$, making $\wtd\scrM$ even ``more invertible''.  In such a situation, determining $E$ by \eqref{eq:SARE-CS''} has an implicit effect of serving the purpose of what an ideal compensated splitting \eqref{eq:t-scrMscrN} should achieve.
This heuristic explanation is actually reflected in our numerical experiments in \cref{sec:egs}.

Now with the optimal $E$ of \eqref{eq:SARE-CS''} in \Cref{thm:GCS(complex)}, we create
a compensated splitting in \eqref{eq:t-scrMscrN} and numerically the fixed-point
iterative scheme \eqref{eq:NatItn-pert} which we will call the {\em compensated fixed-point iteration} (cFPI), as outlined in \Cref{alg:SPLT-global}, where implementation details are not included in order to keep it concise. We now discuss some  implementation details of \Cref{alg:SPLT-global}:
\begin{enumerate}[(1)]
  \item All inner Lyapunov equations at Line 5 share the same coefficient matrix $A+E$ and, hence,
        if solved by the Bartels-Stewart algorithm \cite{bast:1972},
        the Schur decomposition of $A+E$ should be computed outside the for-loop. Further if
        $A+E$ is a real matrix, one should go for its real Schur decomposition.
        As a result of the Schur decomposition, it is rather convenient to check if the operator
        $\wtd\scrM$ defined as in \eqref{eq:NatItn-pert} is invertible or not, namely
        $\wtd\scrM$ is invertible if $\lambda+\bar\mu\ne 0$ for any eigenvalues $\lambda,\,\mu$ of $A+E$.
        cFPI cannot go forward if $\wtd\scrM$ is not invertible. For the purpose of numerical experiments,
        this is what we will do in \cref{sec:egs}.
  \item The normalized residual $\NRes(\cdot)$ at Line 6 is defined as
\begin{equation}\label{eq:residual}
\NRes(X)
   :=\frac {\Big\|A X+ X A^{\HH} + \sum_{j=1}^{m}N_j X N_j^{\HH} + BB^{\HH}\Big\|_{\F}}
           {\|X\|_{\F}(2\|A\|_2+\sum_{j=1}^m\|N_j\|_2^2)+\|BB^{\HH}\|_{\F}}.
\end{equation}
Since the spectral norm $\|\cdot\|_2$ is computationally nontrivial, we may replace
$\|A\|_2$ with $\sqrt{\|A\|_1\|A\|_{\infty}}$, and similarly for $\|N_j\|_2$.
These replacements should work equally well numerically because what is important is the
correct magnitude of $\NRes(\cdot)$.
  \item As we mentioned earlier, we are not directly minimize $\rho(\wtd\scrM^{-1}\wtd\scrN\,)$
        subject to $\wtd\scrM$ being invertible, because there is no easy way to do so.
        Also, without additional efforts, we do not know if $\rho(\wtd\scrM^{-1}\wtd\scrN\,)<1$
        or not. However, there are a couple of remedies. One is
        to estimate $\rho(\wtd\scrM^{-1}\wtd\scrN\,)$ by
        Krylov-type methods \cite{bddrv:2000}, such as MATLAB's {\tt eigs} \cite{lesy:1998}, to roughly calculate the dominate eigenvalue in absolute value of $\wtd\scrM^{-1}\wtd\scrN$. This is especially feasible if
        the Schur decomposition of $A+E$ is available, making the matrix-vector product with $\wtd\scrM^{-1}\wtd\scrN$
        convenient to do. Another one, which is primitive, is to run a certain number of cFPI iterations and then check
        if $\NRes(X_k)$ is on an eventually declining trend. If yes, let the iteration continue until
        the stopping criterion $\NRes(X_k)\le\epsilon$ is satisfied; otherwise, declaim that cFPI fails. It is noted
        that often in actually computations, one usually set a maximal number of iterations to avoid
        a potentially endless loop anyway.
\end{enumerate}

\begin{algorithm}[t]
\caption{cFPI: the compensated fixed-point iteration for \eqref{eq:genlyap}}
\label{alg:SPLT-global}
\begin{algorithmic}[1]
\REQUIRE $A$, $N_j$, $B$ as in \eqref{eq:genlyap}, $\mu\ge 0$, and tolerance $\epsilon$.
\ENSURE  an approximation to the solution of \eqref{eq:genlyap}.
\STATE compute $E$ according to \eqref{eq:optimal(E):cpx};
\STATE $X_0 = 0$.
\FOR{$k = 1,2,\dots$}
\STATE $C=-\sum\limits_{j=1}^{m}N_j X_{k-1} N_j^{\HH} + EX_{k-1} + X_{k-1}E^{\HH} - BB^{\HH}$;
\STATE solve  $(A+E)X_k + X_k (A+E)^{\HH} =  C$  for $X_k$;
\IF{ $\NRes(X_k)\le\epsilon$ }
   \RETURN $X_k$.
\ENDIF
\ENDFOR
\end{algorithmic}
\end{algorithm}

In the next two subsections, we will establish \eqref{eq:optimal(E):cpx}.
In theory, the real case is a special one of the complex case when all involved matrices are real matrices.
Hence there is no need to consider the real case separately if we can handle the complex case. It turns out
that the real case is so much simpler and often the generalized Lyapunov equations \eqref{eq:genlyap} in applications
are real ones, so we will single out the real case first and dive into the complex case separately.
Henceforth, $E_{(s,t)}$ stands for the $(s,t)$th entry of $E$ and similarly for $(N_j)_{(s,t)}$.

\subsection{Optimal $E$ of \eqref{eq:SARE-CS''} for the Real Case}\label{ssec:E-real-proof}

Consider $E\in\bbR^{n\times n}$ and so are all $N_j$. We have
\begin{align}\label{eq:obj:real}
f(E)&=\sum_{s=1}^n\Bigl\|E+E_{(s,s)} I_n-\sum_{j=1}^m(N_j)_{(s,s)}N_j\Bigr\|_{\F}^2 \nonumber\\
    &\quad +\sum_{\begin{subarray}{c}
                  s,t=1\\ s\ne t
                 \end{subarray}}^n\Bigl\|E_{(s,t)} I_n-\sum_{j=1}^m(N_j)_{(s,t)} N_j\Bigr\|_{\F}^2
           +\mu\sum_{s,t=1}^n E_{(s,t)}^2.
\end{align}
We start by finding all possible stationary points of \eqref{eq:SARE-CS''}. To that end, it suffices to set the partial derivatives of $f$ with respect to each entry of $E$ to $0$.
We have,
for $s,t=1,2,\ldots, n$ and $s\ne t$,
\begin{subequations}\label{eq:partialf(E):real}
\begin{align}
0&=\frac 14\frac {\partial f(E)}{\partial E_{(s,s)}}\nonumber\\
  &=\frac 14\frac {\partial }{\partial E_{(s,s)}}\Bigl\|E+E_{(s,s)} I_n-\sum_{j=1}^m(N_j)_{(s,s)}N_j\Bigr\|_{\F}^2\nonumber\\
  &\quad+\frac 14\frac {\partial }{\partial E_{(s,s)}}
      \sum_{i=1,\, i\ne s}^n\Bigl\|E+E_{(i,i)} I_n-\sum_{j=1}^m(N_j)_{(i,i)}N_j\Bigr\|_{\F}^2
      +\frac{\mu}{4}\frac{\partial}{\partial E_{(s,s)}} E_{(s,s)}^2\nonumber\\
  &=\frac 12\Bigl[2E_{(s,s)}-\sum_{j=1}^m\big[(N_j)_{(s,s)}\big]^2\Bigr]2\nonumber\\
  &\quad+\frac 12\sum_{i=1,\,i\ne s}^n
       \Bigl[E_{(i,i)}+E_{(s,s)}-\sum_{j=1}^m(N_j)_{(s,s)}(N_j)_{(i,i)}\Bigr]\nonumber\\
  &\quad+\frac 12\sum_{i=1,\,i\ne s}^n
    \Bigl[E_{(s,s)}+E_{(i,i)}-\sum_{j=1}^m(N_j)_{(i,i)}(N_j)_{(s,s)}\Bigr]
       +\frac{\mu}{2}E_{(s,s)}\nonumber\\
  &=\left(n+1+\frac{\mu}{2}\right)E_{(s,s)}
       +\sum_{i=1,\,i\ne s}^n E_{(i,i)}
       -\sum_{j=1}^m(N_j)_{(s,s)}\sum_{i=1}^n(N_j)_{(i,i)}\nonumber\\
  &=\left(n+\frac{\mu}{2}\right)E_{(s,s)}
       +\sum_{i=1}^n E_{(i,i)}
       -\sum_{j=1}^m(N_j)_{(s,s)}\tr(N_j),
       \label{eq:partialf(E):diag:real}\\
0&=\frac 12\frac {\partial f(E)}{\partial E_{(s,t)}}\nonumber\\
  &=\frac 12\frac {\partial }{\partial E_{(s,t)}}
          \sum_{i=1}^n\Bigl\|E+E_{(i,i)} I_n-\sum_{j=1}^m(N_j)_{(i,i)}N_j\Bigr\|_{\F}^2\nonumber\\
  &\quad+\frac 12\frac {\partial }{\partial E_{(s,t)}}
          \Bigl\|E_{(s,t)} I_n-\sum_{j=1}^m(N_j)_{(s,t)}N_j\Bigr\|_{\F}^2
          +\frac{\mu}{2}\frac{\partial }{\partial E_{(s,t)}}E_{(s,t)}^2\nonumber\\
  &=\sum_{i=1}^n\Bigl[E_{(s,t)}-\sum_{j=1}^m(N_j)_{(i,i)}(N_j)_{(s,t)}\Bigr]\nonumber\\
  &\quad+\sum_{i=1}^n\Bigl[E_{(s,t)}-\sum_{j=1}^m(N_j)_{(s,t)}(N_j)_{(i,i)}\Bigr]
          +\mu E_{(s,t)}\nonumber\\
  &=(2n+\mu)E_{(s,t)}-2\sum_{j=1}^m(N_j)_{(s,t)}\tr(N_j).
    \label{eq:partialf(E):offdiag:real}
\end{align}
\end{subequations}
First we solve the system of $n$ equations in \eqref{eq:partialf(E):diag:real} for $1\le s\le n$
for the diagonal entries of $E$. Sum up the $n$ equations over
$s=1,2,\ldots, n$ to get
\begin{equation}\label{eq:partialf(E):diag'}
\left(2n+\frac{\mu}{2}\right)\sum_{s=1}^n E_{(s,s)}
-\sum_{j=1}^m[\tr(N_j)]^2=0
\quad\Rightarrow\quad
\sum_{s=1}^n E_{(s,s)}=\frac{2}{4n+\mu}\sum_{j=1}^m[\tr(N_j)]^2.
\end{equation}
Plugging \eqref{eq:partialf(E):diag'} into \eqref{eq:partialf(E):diag:real}, we find
\begin{equation}\label{eq:optimal(E):diag:real}
E_{(s,s)}=\frac{2}{2n+\mu}
        \left[\sum_{j=1}^m(N_j)_{(s,s)}\tr(N_j)
              -\frac{2}{4n+\mu}\sum_{j=1}^m[\tr(N_j)]^2\right]\quad
        \mbox{for $s=1,2,\ldots, n$.}
\end{equation}
Consider now the off-diagonal entries of $E$. By \eqref{eq:partialf(E):offdiag:real}, we have
for $s,\,t=1,2,\ldots, n$ and $s\ne t$, yielding
\begin{equation}\label{eq:optimal(E):offdiag:real}
E_{(s,t)}=\frac{2}{2n+\mu}\sum_{j=1}^m(N_j)_{(s,t)}\tr(N_j)\quad
        \mbox{for $s,\,t=1,2,\ldots, n$ and $s\ne t$.}
\end{equation}
Putting \eqref{eq:optimal(E):diag:real} and \eqref{eq:optimal(E):offdiag:real} together, compactly, we have
\eqref{eq:optimal(E):cpx} for the real case.
For every $\mu\ge 0$, there is exactly one stationary point, given by
\eqref{eq:optimal(E):diag:real} and \eqref{eq:optimal(E):offdiag:real}.
By convexity, it is the unique global minimizer.
%
%
%

\subsection{Optimal $E$ of \eqref{eq:SARE-CS''} for the Complex Case}\label{ssec:E-complex-proof}
Consider the case when $E$ and all $N_j$ are complex.
We have
\begin{align*}
f(E)
  &=\sum_{s=1}^n\Bigl\|E+\bE_{(s,s)} I_n-\sum_{j=1}^m(\bN_j)_{(s,s)}N_j\Bigr\|_{\F}^2\\
  &\quad +\sum_{\begin{subarray}{c}
                s,t=1\\ s\ne t
                \end{subarray}}^n\Bigl\|\bE_{(s,t)} I_n-\sum_{j=1}^m(\bN_j)_{(s,t)} N_j\Bigr\|_{\F}^2
    +\mu\sum_{s,t=1}^n E_{(s,t)}\bE_{(s,t)}.
\end{align*}
To find the optimal $E$ to \eqref{eq:SARE-CS''}, we will regard $f(E)$ as a function of the real and imaginary parts of
all entries of $E$ and find all possible stationary points of \eqref{eq:SARE-CS''} determined by setting
the partial derivatives of $f$ with respect to the real and imaginary parts of each entry of $E$ to $0$.
But first we will need to derive those partial derivatives.
To that end, it is rather convenient to do so via through Wirtinger derivatives \cite{koor2023shorttutorialwirtingercalculus} by
introducing
$$
F(E,\bar E):=f(E)
$$
and at the same time treating $F(E,\bar E)$
as a function of two independent matrix variables $E$ and $\bar E$ (although $\bar E$ is the entrywise  complex conjugate of $E$).
%
%
Write
$$
E_{(s,t)} = E_{(s,t)}^{{\rm re}} + \iota E_{(s,t)}^{{\rm im}}, \quad
\bE_{(s,t)} = E_{(s,t)}^{{\rm re}} - \iota E_{(s,t)}^{{\rm im}},
$$
where $\iota$ is the imaginary unit, and $E_{(s,t)}^{{\rm re}},\,E_{(s,t)}^{{\rm im}}\in\bbR$.
It is known that
\begin{equation}\label{eq:diff-realE-imagE}
\frac {\partial f(E)}{\partial E_{(s,t)}^{{\rm re}}}
= \frac {\partial F(E,\bar E)}{\partial E_{(s,t)}} + \frac {\partial F(E,\bar E)}{\partial \bE_{(s,t)}}, 
\quad
\frac 1{\iota}\cdot\frac {\partial f(E)}{\partial E_{(s,t)}^{{\rm im}}}
= \frac {\partial F(E,\bar E)}{\partial E_{(s,t)}} - \frac {\partial F(E,\bar E)}{\partial \bE_{(s,t)}}. 
\end{equation}
For the diagonal entries, we have
\begin{subequations}\label{eq:partial-diag:cpx}
\begin{align}
\frac {\partial F(E,\bar E)}{\partial E_{(s,s)}}
  &=\frac {\partial }{\partial E_{(s,s)}}
        \Bigl\|E+\bE_{(s,s)} I_n-\sum_{j=1}^m(\bN_j)_{(s,s)}N_j\Bigr\|_{\F}^2 \nonumber \\
  &\quad   + \frac {\partial }{\partial E_{(s,s)}}
       \sum_{i=1,\,i\ne s}^n
           \Bigl\|E+\bE_{(i,i)} I_n-\sum_{j=1}^m(\bN_j)_{(i,i)}N_j\Bigr\|_{\F}^2
       +\mu\frac {\partial (E_{(s,s)}\bE_{(s,s)})}{\partial E_{(s,s)}} \nonumber\\
  &= \frac {\partial }{\partial E_{(s,s)}}
       \Bigl| E_{(s,s)} + \bE_{(s,s)} - \sum_{j=1}^m (\bN_j)_{(s, s)}(N_j)_{(s, s)}\Bigr|^2\nonumber \\
  &\quad+\frac {\partial }{\partial E_{(s,s)}}
       \sum_{i=1,\,i\ne s}^n
        \Bigl|E_{(i, i)} + \bE_{(s,s)} - \sum_{j=1}^m (\bN_j)_{(s, s)}(N_j)_{(i,i)} \Bigr|^2 \nonumber \\
  &\quad+\mu\frac {\partial }{\partial E_{(s,s)}}
       \Bigl| E_{(s,s)} + \bE_{(s,s)}\Bigr|^2
         + \mu\frac {\partial }{\partial E_{(s,s)}}
                  \sum_{i=1,\,i\ne s}^n
               \Bigl|E_{(i, i)} + \bE_{(s,s)} \Bigr|^2 \nonumber \\
    &\quad  + \frac {\partial }{\partial E_{(s,s)}}
                \sum_{i=1,\,i\ne s}^n
       \Bigl|E_{(s,s)}+\bE_{(i,i)}-\sum_{j=1}^m(\bN_j)_{(i,i)}(N_j)_{(s,s)}\Bigr|^2
       +\mu\bE_{(s,s)} \nonumber\\
  &=2\Bigl[E_{(s,s)}+\bE_{(s,s)}
       -\sum_{j=1}^m(\bN_j)_{(s,s)}(N_j)_{(s,s)}\Bigr] \nonumber\\
  &\quad+\sum_{i=1,\,i\ne s}^n
       \Bigl[E_{(i,i)}+\bE_{(s,s)}-\sum_{j=1}^m(\bN_j)_{(s,s)}(N_j)_{(i,i)}\Bigr] \nonumber\\
  &\quad+\sum_{i=1,\,i\ne s}^n
       \Bigl[\bE_{(s,s)}+E_{(i,i)}-\sum_{j=1}^m(N_j)_{(i,i)}(\bN_j)_{(s,s)}\Bigr]
       +\mu\bE_{(s,s)} \nonumber\\
  &=2\sum_{i=1}^n
       \Bigl[\bE_{(s,s)}+E_{(i,i)}-\sum_{j=1}^m(N_j)_{(i,i)}(\bN_j)_{(s,s)}\Bigr]
       +\mu\bE_{(s,s)} \nonumber\\
  &=(2n+\mu)\bE_{(s,s)}+2\sum_{i=1}^nE_{(i,i)}
      -2\sum_{j=1}^m \tr(N_j)\,(\bN_j)_{(s, s)}   ,     \label{eq:partial-diag:1:cpx} \\
\frac {\partial F(E,\bar E)}{\partial \bE_{(s,s)}}
  & =\frac {\partial }{\partial \bE_{(s,s)}}
        \Bigl\|E+\bE_{(s,s)} I_n-\sum_{j=1}^m(\bN_j)_{(s,s)}N_j\Bigr\|_{\F}^2 \nonumber\\
  &\quad+\frac {\partial }{\partial \bE_{(s,s)}}
       \sum_{i=1,\,i\ne s}^n
           \Bigl\|E+\bE_{(i,i)} I_n-\sum_{j=1}^m(\bN_j)_{(i,i)}N_j\Bigr\|_{\F}^2
       +\mu\frac {\partial (E_{(s,s)}\bE_{(s,s)})}{\partial \bE_{(s,s)}} \nonumber\\
  & =\frac {\partial }{\partial \bE_{(s,s)}}
     \Bigl|E_{(s, s)} + \bE_{(s, s)} - \sum_{j=1}^m (\bN_j)_{(s, s)}(N_j)_{(s, s)}\Bigr|^2 \nonumber \\
  &\quad+\frac {\partial }{\partial \bE_{(s,s)}}
       \sum_{i=1,\,i\ne s}^n
    \Bigl|E_{(i, i)} + \bE_{(s,s)} - \sum_{j=1}^m (\bN_j)_{(s,s)}(N_j)_{(i,i)} \Bigr|^2\nonumber \\
  &\quad+\frac {\partial }{\partial \bE_{(s,s)}}
       \sum_{i=1,\,i\ne s}^n
       \Bigl|\bE_{(s,s)}+E_{(i,i)}-\sum_{j=1}^m(N_j)_{(i,i)}(\bN_j)_{(s,s)}\Bigr|^2
       +\mu E_{(s,s)} \nonumber\\
  &=2\Bigl[E_{(s,s)}+\bE_{(s,s)}
       -\sum_{j=1}^m(\bN_j)_{(s,s)}(N_j)_{(s,s)}\Bigr] \nonumber\\
  &\quad+\sum_{i=1,\,i\ne s}^n
       \Bigl[\bE_{(i,i)}+E_{(s,s)}-\sum_{j=1}^m(N_j)_{(s,s)}(\bN_j)_{(i,i)}\Bigr] \nonumber\\
  &\quad+\sum_{i=1,\,i\ne s}^n
       \Bigl[E_{(s,s)}+\bE_{(i,i)}-\sum_{j=1}^m(\bN_j)_{(i,i)}(N_j)_{(s,s)}\Bigr]
       +\mu E_{(s,s)} \nonumber\\
  &=2\sum_{i=1}^n
       \Bigl[E_{(s,s)}+\bE_{(i,i)}-\sum_{j=1}^m(\bN_j)_{(i,i)}(N_j)_{(s,s)}\Bigr]
       +\mu E_{(s,s)} \nonumber\\
  &=(2n+\mu)E_{(s,s)}+2\sum_{i=1}^n\bE_{(i,i)}
      -2\sum_{j=1}^m \tr(\bN_j)\,(N_j)_{(s, s)}. \label{eq:partial-diag:2:cpx}
\end{align}
\end{subequations}
According to \eqref{eq:diff-realE-imagE} and upon using \eqref{eq:partial-diag:cpx}, the first order optimality condition on the diagonal entries of $E$
is given by, for $1\le s\le n$,
\begin{subequations}\label{eq:KKT-diag-cpx}
\begin{align}
0&=\frac {\partial f(E)}{\partial E_{(s,s)}^{{\rm re}}}
   =\frac {\partial F(E,\bar E)}{\partial E_{(s,s)}}+\frac {\partial F(E,\bar E)}{\partial \bE_{(s,s)}} \nonumber\\
  &=(2n+\mu)\big[E_{(s,s)}+\bE_{(s,s)}\big]
       +2\sum_{i=1}^n\big[E_{(i,i)}+\bE_{(i,i)}\big] \nonumber\\
  &\quad-4\sum_{j=1}^m\Re\big(\tr(\bN_j)\,(N_j)_{(s,s)}\big) \nonumber\\
  &=(4n+2\mu)E_{(s,s)}^{{\rm re}}+4\sum_{i=1}^nE_{(i,i)}^{{\rm re}}
      -4\sum_{j=1}^m\Re\big( \tr(\bN_j)\,(N_j)_{(s, s)}\big), \label{eq:KKT-diag-Re} \\
0&=\frac 1{\iota}\cdot\frac {\partial f(E)}{\partial E_{(s,s)}^{{\rm im}}}
   =\frac {\partial F(E,\bar E)}{\partial E_{(s,s)}}-\frac {\partial F(E,\bar E)}{\partial \bE_{(s,s)}} \nonumber\\
  &=-(2n+\mu)\big[E_{(s,s)}-\bE_{(s,s)}\big]
       +2\sum_{i=1}^n\big[E_{(i,i)}-\bE_{(i,i)}\big] \nonumber\\
  &\quad+4\iota\sum_{j=1}^m\Im\big(\tr(\bN_j)\,(N_j)_{(s,s)}\big) \nonumber\\
  &=-(4n+2\mu)\iota E_{(s,s)}^{{\rm im}}+4\iota\sum_{i=1}^nE_{(i,i)}^{{\rm im}}
    +4\iota\sum_{j=1}^m\Im\big(\tr(\bN_j)\,(N_j)_{(s, s)}\big). \label{eq:KKT-diag-Im}
\end{align}
\end{subequations}

First we solve \eqref{eq:KKT-diag-Re} for $1\le s\le n$. Summing up \eqref{eq:KKT-diag-Re} over $1\le s\le n$,
we get, after simplifications,
\begin{align}
&(8n+2\mu)\sum_{i=1}^n E_{(i,i)}^{{\rm re}}
   = 4\sum_{j=1}^m\tr(\bN_j) \tr(N_{j}) \nonumber \\
\Rightarrow\quad&
\sum_{i=1}^n E_{(i,i)}^{{\rm re}}
   = \frac{2}{4n+\mu}\sum_{j=1}^m \tr(\bN_j) \tr(N_{j}).
   \label{eq:KKT-diag-Re:sum}
\end{align}
Plug \eqref{eq:KKT-diag-Re:sum} into \eqref{eq:KKT-diag-Re} to get, for $s=1,2,\ldots,n$,
\begin{equation}\label{eq:optimal(E):diag-Re}
E_{(s,s)}^{{\rm re}}
 =\frac{2}{2n+\mu}\left[\sum_{j=1}^m\Re\big((N_j)_{(s,s)}\tr(\bN_j)\big)
 -\frac{2}{4n+\mu}\sum_{j=1}^m\tr(\bN_j)\tr(N_j)\right].
\end{equation}
Next, we solve \eqref{eq:KKT-diag-Im} for $1\le s\le n$, which can be simplified to
\begin{equation}\label{eq:KKT-diag-Im'}
-\Bigl(n+\frac{\mu}{2}\Bigr)E_{(s,s)}^{{\rm im}}+\sum_{i=1}^nE_{(i,i)}^{{\rm im}}
 =-\sum_{j=1}^m\Im\big(\tr(\bN_j)\,(N_j)_{(s,s)}\big),
 \quad 1\le s\le n,
\end{equation}
or, compactly,
\begingroup
\setlength{\arraycolsep}{3.5pt}
\begin{equation}\label{eq:KKT-diag-Im''}
\begin{bmatrix}
n-1+\mu/2 & -1 & \cdots & -1\\
-1 & n-1+\mu/2 & \cdots & -1\\
\vdots          &   \vdots   & \ddots & \vdots \\
-1 & -1 & \cdots & n-1+\mu/2
\end{bmatrix}
\begin{bmatrix}
E_{(1,1)}^{{\rm im}} \\
E_{(2,2)}^{{\rm im}} \\
\vdots \\
E_{(n,n)}^{{\rm im}}
\end{bmatrix}
=\sum_{j=1}^m
\begin{bmatrix}
\Im\big((N_j)_{(1,1)}\tr(\bN_j)\big)\\
\Im\big((N_j)_{(2,2)}\tr(\bN_j)\big)\\
    \vdots \\
\Im\big((N_j)_{(n,n)}\tr(\bN_j)\big)
\end{bmatrix}:=\bb.
\end{equation}
\endgroup
It coefficient matrix is
$$
C=\left(n+\frac {\mu}2\right) I_n-\be\be^{\T}=\frac {2n+\mu}2\left(I_n-\frac 2{2n+\mu}\be\be^{\T}\right)
$$
which is positive definite for $\mu>0$ and positive semidefinite for $\mu=0$,
where $\be\in\bbR^n$ is the vector of all ones.
To solve \eqref{eq:KKT-diag-Im''}, we distinguish the following two cases.

\smallskip
\noindent\textbf{Case 1: $\mu>0$.}
Now $C$ is invertible
and\footnote {It is well-known that $I_n-\bx\by^{\T}$ is invertible if and only if $\by^{\T}\bx\ne 1$, in which case
    $$
    (I_n-\bx\by^{\T})^{-1}=I_n+\frac 1{1-\by^{\T}\bx}\,\bx\by^{\T}.
    $$}
$$
C^{-1}=\frac 2{2n+\mu}\left(I_n+\frac 2{\mu}\be\be^{\T}\right).
$$
Noticing also that
\begin{align}
\be^{\T}\bb
 &=\sum_{i=1}^n\sum_{j=1}^m
   \Im\big((N_j)_{(i,i)}\tr(\bN_j)\big) \nonumber \\
 &=\sum_{j=1}^m
   \Im\big(\tr(N_j)\tr(\bN_j)\big) \nonumber\\
 &=\sum_{j=1}^m\Im\big(|\tr(N_j)|^2\big)=0, \label{eq:eTb=0}
\end{align}
we find the unique solution to \eqref{eq:KKT-diag-Im''} is thus $\frac 2{2n+\mu}\bb$, i.e.,
\begin{equation}\label{eq:optimal(E):diag-Im:positive}
E_{(s,s)}^{{\rm im}}
 =\frac{2}{2n+\mu}\sum_{j=1}^m
 \Im\big((N_j)_{(s,s)}\tr(\bN_j)\big)
 \quad\mbox{for $1\le s\le n$}.
\end{equation}

\smallskip
\noindent\textbf{Case 2: $\mu=0$.}
The coefficient matrix $C$ of \eqref{eq:KKT-diag-Im''} is singular.
We make a couple of observations:
\begin{enumerate}[1)]
  \item the linear system \eqref{eq:KKT-diag-Im''} is consistent because
the range of $C$ is the orthogonal complement of ${\rm span}\{\be\}$ and $\bb$ is in the orthogonal complement
because of \eqref{eq:eTb=0};
  \item Since that $I_n-(1/n)\be\be^{\T}$ is an orthogonal projection matrix, its Moore-Penrose
pseudo-inverse \cite[p.290]{govl:2013} is itself and hence the Moore-Penrose pseudo-inverse of $C$ is given by
$C^{\dagger}=\frac 1n\Big(I_n-\frac 1n\be\be^{\T}\Big)$, yielding the minimum norm solution to \eqref{eq:KKT-diag-Im''}:
        $
        C^{\dagger}\bb=\frac 1n\bb,
        $
        i.e.,
$$ 
E_{(s,s)}^{{\rm im}}
 =\frac 1n\sum_{j=1}^m\Im\big((N_j)_{(s,s)} \tr(\bN_{j})\big).
$$ 
\end{enumerate}
Consequently, the unique solution for $\mu>0$ and the minimum norm solution for $\mu=0$
can be expressed as
\begin{equation}\label{eq:optimal(E):diag-Im}
E_{(s,s)}^{{\rm im}}
 =\frac{2}{2n+\mu}\sum_{j=1}^m\Im\big((N_j)_{(s,s)}\tr(\bN_j)\big)
 \quad\mbox{for $1\le s\le n$}.
\end{equation}
Putting \eqref{eq:optimal(E):diag-Re} and \eqref{eq:optimal(E):diag-Im} together gives
\begin{equation}\label{eq:optimal(E):diag:cpx}
E_{(s,s)}=\frac{2}{2n+\mu}\left[\sum_{j=1}^m(N_j)_{(s,s)}\tr(\bN_j)
 -\frac{2}{4n+\mu}\sum_{j=1}^m\tr(\bN_j)\tr(N_j)\right].
\end{equation}

Now, we consider the off-diagonal entries of $E$ and $\bE$.
For $s\ne t$, we have
\begin{subequations}\label{eq:partial-offdiag:cpx}
\begin{align}
\frac {\partial F(E,\bar E)}{\partial E_{(s,t)}}
  &=\sum_{i=1}^n\frac {\partial }{\partial E_{(s,t)}}
       \Bigl|E_{(s,t)}-\sum_{j=1}^m(\bN_j)_{(i,i)}(N_j)_{(s,t)}\Bigr|^2 \nonumber\\
  &\quad+\sum_{i=1}^n\frac {\partial }{\partial E_{(s,t)}}
             \Bigl|\bE_{(s,t)} -  \sum_{j=1}^m (\bN_j)_{(s,t)}(N_j)_{(i,i)}\Bigr|^2
       +\mu\frac {\partial(E_{(s,t)}\bE_{(s,t)})}{\partial E_{(s,t)}} \nonumber\\
  &=\sum_{i=1}^n\Bigl(\bE_{(s,t)}
       -\sum_{j=1}^m(N_j)_{(i,i)}(\bN_j)_{(s,t)}\Bigr) \nonumber\\
  &\quad+\sum_{i=1}^n\Bigl(\bE_{(s,t)}
       -\sum_{j=1}^m(\bN_j)_{(s,t)}(N_j)_{(i,i)}\Bigr)
       +\mu\bE_{(s,t)} \nonumber\\
  &=(2n+\mu)\bE_{(s,t)}
       -2\sum_{j=1}^m(\bN_j)_{(s,t)}\tr(N_j), \label{eq:partial-offdiag:1:cpx}\\
\frac {\partial F(E,\bar E)}{\partial \bE_{(s,t)}}
   &= \sum_{i=1}^n \frac {\partial }{\partial \bE_{(s,t)}}
       \Bigl|E_{(s,t)}-\sum_{j=1}^m(\bN_j)_{(i,i)}(N_j)_{(s,t)}\Bigr|^2 \nonumber\\
  &\quad+\sum_{i=1}^n\frac {\partial }{\partial \bE_{(s,t)}}
         \Bigl|\bE_{(s,t)} -  \sum_{j=1}^m (\bN_j)_{(s,t)}(N_j)_{(i,i)} \Bigr|^2
       +\mu\frac {\partial(E_{(s,t)}\bE_{(s,t)})}{\partial \bE_{(s,t)}} \nonumber\\
  &=\sum_{i=1}^n\Bigl(E_{(s,t)}
       -\sum_{j=1}^m(\bN_j)_{(i,i)}(N_j)_{(s,t)}\Bigr) \nonumber\\
  &\quad+\sum_{i=1}^n\Bigl(E_{(s,t)}
       -\sum_{j=1}^m(N_j)_{(s,t)}(\bN_j)_{(i,i)}\Bigr)
       +\mu E_{(s,t)} \nonumber\\
  &=(2n+\mu)E_{(s,t)}
       -2\sum_{j=1}^m(N_j)_{(s,t)}\tr(\bN_j). \label{eq:partial-offdiag:2:cpx}
\end{align}
\end{subequations}
According to \eqref{eq:diff-realE-imagE}  and upon using \eqref{eq:partial-offdiag:cpx}, the first order optimality
condition on the off-diagonal entries of $E$ is given by,
for $1\le s,\,t\le n$ and $s\ne t$,
\begin{subequations}\label{eq:KKT-offdiag-cpx}
\begin{align}
0&=\frac {\partial f(E)}{\partial E_{(s,t)}^{{\rm re}}}
   = \frac {\partial F(E,\bar E)}{\partial E_{(s,t)}} + \frac {\partial F(E,\bar E)}{\partial \bE_{(s,t)}} \nonumber\\
   &= (2n+\mu) \big[\bE_{(s,t)} + E_{(s,t)}\big] \nonumber\\
   &\quad - 2\sum_{j=1}^m\Bigl(  (\bN_j)_{(s,t)}\tr(N_j) +  (N_j)_{(s,t)}\tr(\bN_j)\Bigr)
         \nonumber\\
   &=(4n+2\mu)E_{(s,t)}^{{\rm re}}-4\sum_{j=1}^m\Re\big((\bN_j)_{(s,t)}\tr(N_j)\big)  , \label{eq:KKT-offdiag-Re}\\
0&=\frac 1{\iota}\cdot\frac {\partial f(E)}{\partial E_{(s,t)}^{{\rm im}}}
   =  \frac {\partial F(E,\bar E)}{\partial E_{(s,t)}} - \frac {\partial F(E,\bar E)}{\partial \bE_{(s,t)}} \nonumber\\
   &=(2n+\mu) \big[\bE_{(s,t)} - E_{(s,t)}\big] \nonumber\\
   &\quad - 2\sum_{j=1}^m \Bigl((\bN_j)_{(s,t)}\tr(N_j)  - (N_j)_{(s,t)}\tr(\bN_j) \Bigr)\nonumber\\
   &=-(4n+2\mu)\iota E_{(s,t)}^{{\rm im}}-4\iota \sum_{j=1}^m\Im\big((\bN_j)_{(s,t)}\tr(N_j)\big). \label{eq:KKT-offdiag-Im}
\end{align}
\end{subequations}
Equations in \eqref{eq:KKT-offdiag-cpx} can be readily solved to yield
$$ 
E_{(s,t)}^{{\rm re}}=\frac{2}{2n+\mu}\sum_{j=1}^m\Re\big((\bN_j)_{(s,t)}\tr(N_j)\big), \quad
E_{(s,t)}^{{\rm im}}
 =-\frac{2}{2n+\mu}\sum_{j=1}^m\Im\big((\bN_j)_{(s,t)}\tr(N_j)\big).
$$ 
Putting them together gives
\begin{equation}\label{eq:optimal(E):offdiag:cpx}
E_{(s,t)}=\frac{2}{2n+\mu}\sum_{j=1}^m(N_j)_{(s,t)}\tr(\bN_j).
\end{equation}

In summary, the real parts $E_{(s,s)}^{{\rm re}}$ of the diagonal entries and all off-diagonal entries $E_{(s,t)}$ for $s\ne t$
are uniquely determined by the first order optimality condition for any $\mu\ge0$, but
the imaginary parts $E_{(s,s)}^{{\rm im}}$ of the diagonal entries are uniquely determined only if $\mu>0$. Hence
there is only one stationary point if $\mu>0$ and that stationary point is the global minimizer.
If, however, $\mu=0$, the first order optimality condition \eqref{eq:KKT-diag-Im''}
for the imaginary parts $E_{(s,s)}^{{\rm im}}$ of the diagonal entries
has infinitely many solutions, and the one with the smallest Frobenius norm is given by
\eqref{eq:optimal(E):diag-Im}.
Finally, putting \eqref{eq:optimal(E):diag:cpx} and \eqref{eq:optimal(E):offdiag:cpx} together,
we get \eqref{eq:optimal(E):cpx}.

\section{Numerical Experiments}\label{sec:egs}

In this section, we will test our proposed method cFPI in \Cref{alg:SPLT-global}, along the way we will compare it against the baseline method of similar complexity: sFPI (the standard fixed point iteration) \eqref{eq:NatItn} from the natural splitting \eqref{eq:NatSplit}.

There will be three experiments: two for real \eqref{eq:genlyap} and one complex problem.
When randomness is involved, we first generate mentioned random matrices and save them for repetitions.
In calling our cFPI (\Cref{alg:SPLT-global}) and sFPI, we stop the iteration as soon as
$\NRes(X_k)\le 10^{-10}$ where $\NRes(\cdot)$ is as
defined in \eqref{eq:residual}. The spectral radii
of $\scrM^{-1}\scrN$ and $\wtd\scrM^{-1}\wtd\scrN$ are estimated by MATLAB's {\tt eigs} \cite{lesy:1998}
upon passing operator-vector product functions by $\scrM^{-1}\scrN$
and $\wtd\scrM^{-1}\wtd\scrN$.

As far as $\rho(\scrM^{-1}\scrN)$ is concerned, it is invariant with respect to substitution $A\leftarrow -A$.
Earlier in \cref{sec:intro}, we commented on that an ideally compensated splitting should be the one that minimizes $\rho(\wtd\scrM^{-1}\wtd\scrN\,)$ subject to $\wtd\scrM$ being invertible, but that is too ambitious to be practical.
We therefore settle on a more practical approach: determine $E$ by \eqref{eq:SARE-CS''} which yields the optimal $E$
as in \eqref{eq:optimal(E):cpx}. This optimal $E$ is also invariant with respect to substitution $A\leftarrow -A$.
Earlier immediately after \Cref{thm:GCS(complex)}, we commented on a potential impact by the eigenvalue
distribution being in the left-half or right-half complex plane.
In order to reflect such an impact,
in each of the following three examples, we actually experiment on two generalized Lyapunov equations \eqref{eq:genlyap}:
one with $A$ and the other with the opposite of it.

At this point, we do not have a way to choose a suitable $\mu$, except letting $\mu$ run from $0$ to some number bigger than $0$.
During testing, we played with $\mu\in [0,0.5]$ on the three experiments. We did not see (much)  smaller $\rho(\wtd\scrM^{-1}\wtd\scrN\,)$
for $\mu>0$ than for $\mu=0$. For that reason, in our following experiments, we take $\mu=0$. How to pick a better $\mu$ will
be left for future research in the days ahead.

\begin{figure}[t]
{\centering
{\resizebox*{0.32\textwidth}{0.17\textheight}{\includegraphics{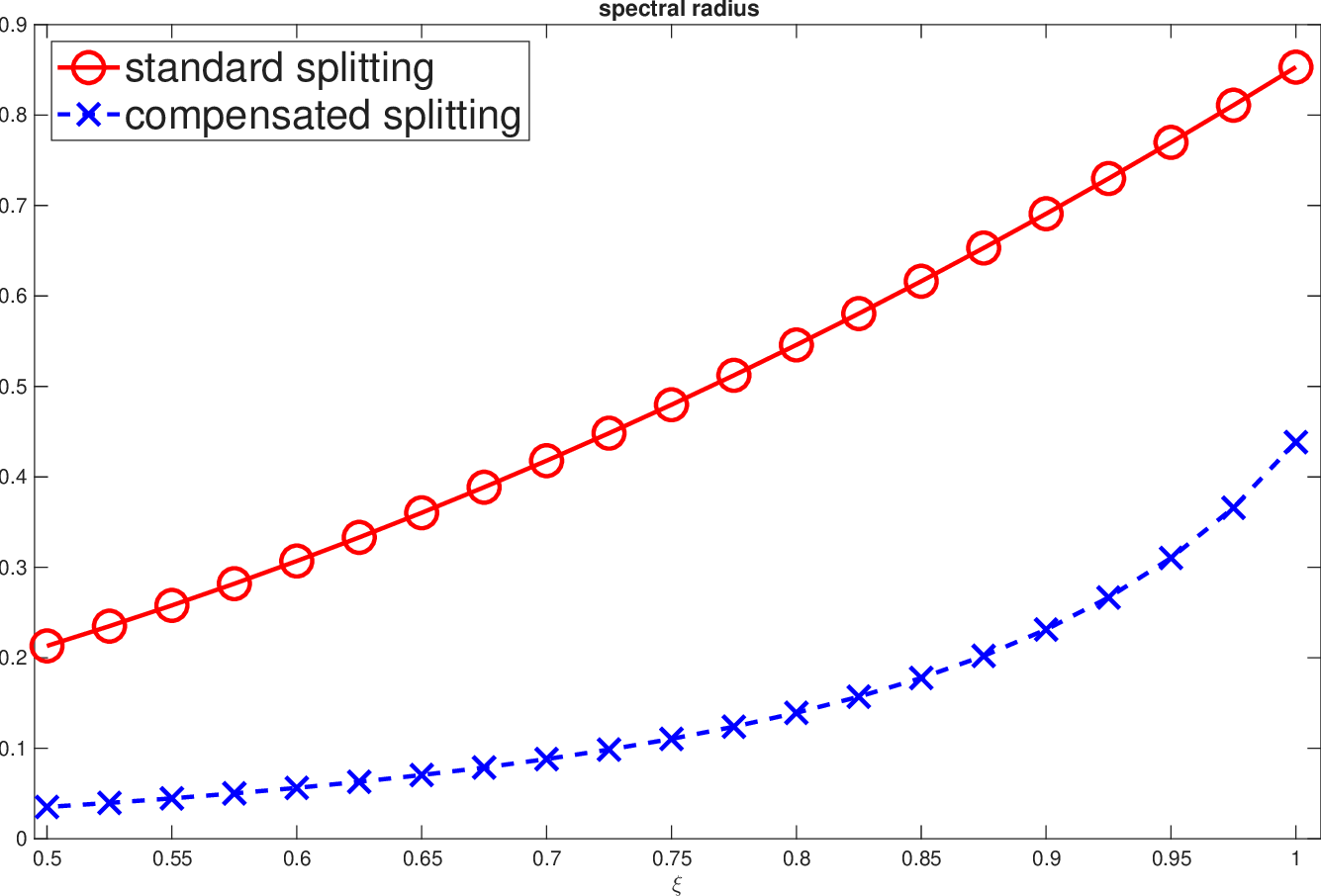}}}\,
{\resizebox*{0.32\textwidth}{0.17\textheight}{\includegraphics{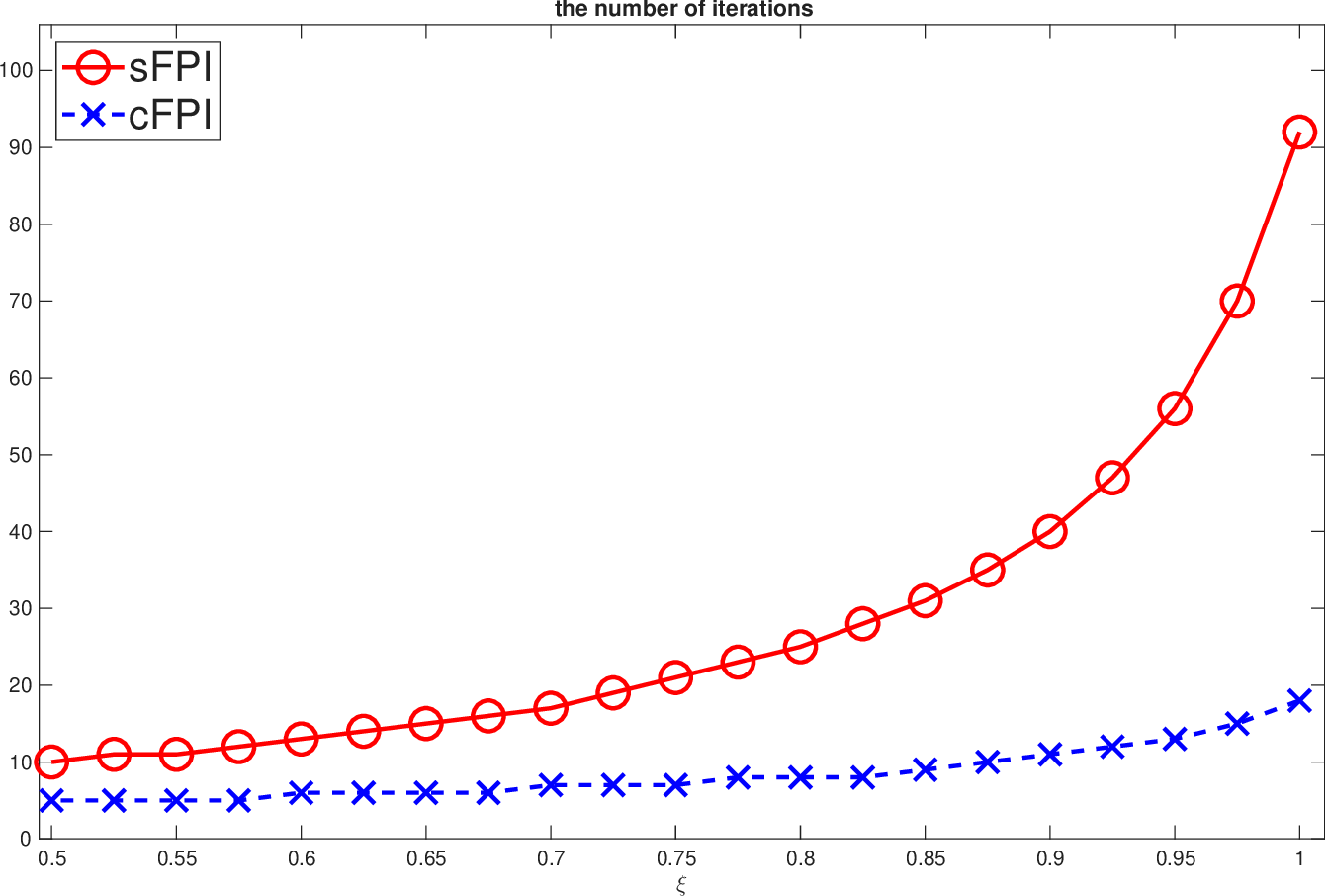}}}\,
{\resizebox*{0.32\textwidth}{0.17\textheight}{\includegraphics{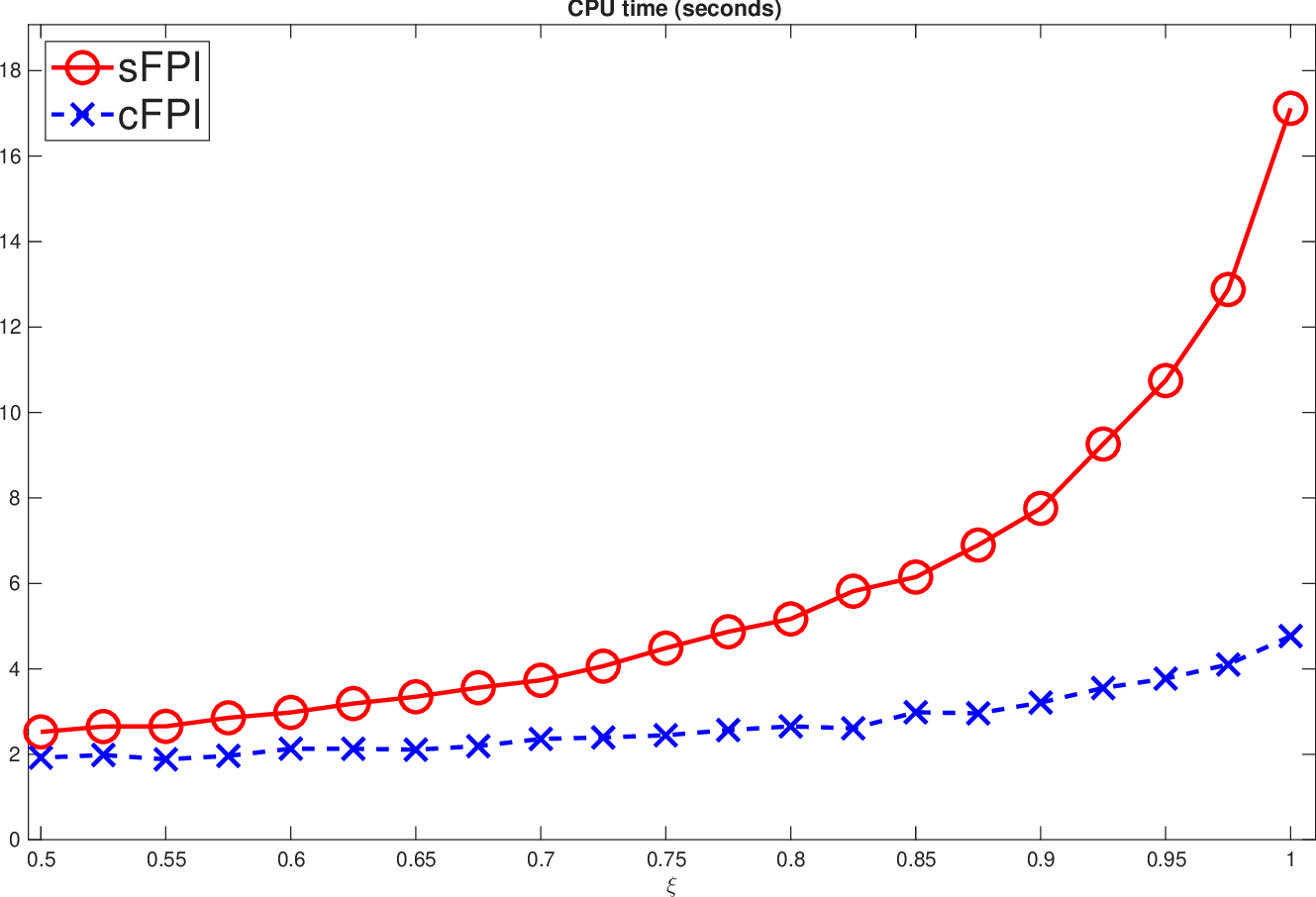}}}
} \\
{\centering
\subfloat[spectral radius]{\resizebox*{0.32\textwidth}{0.17\textheight}{\includegraphics{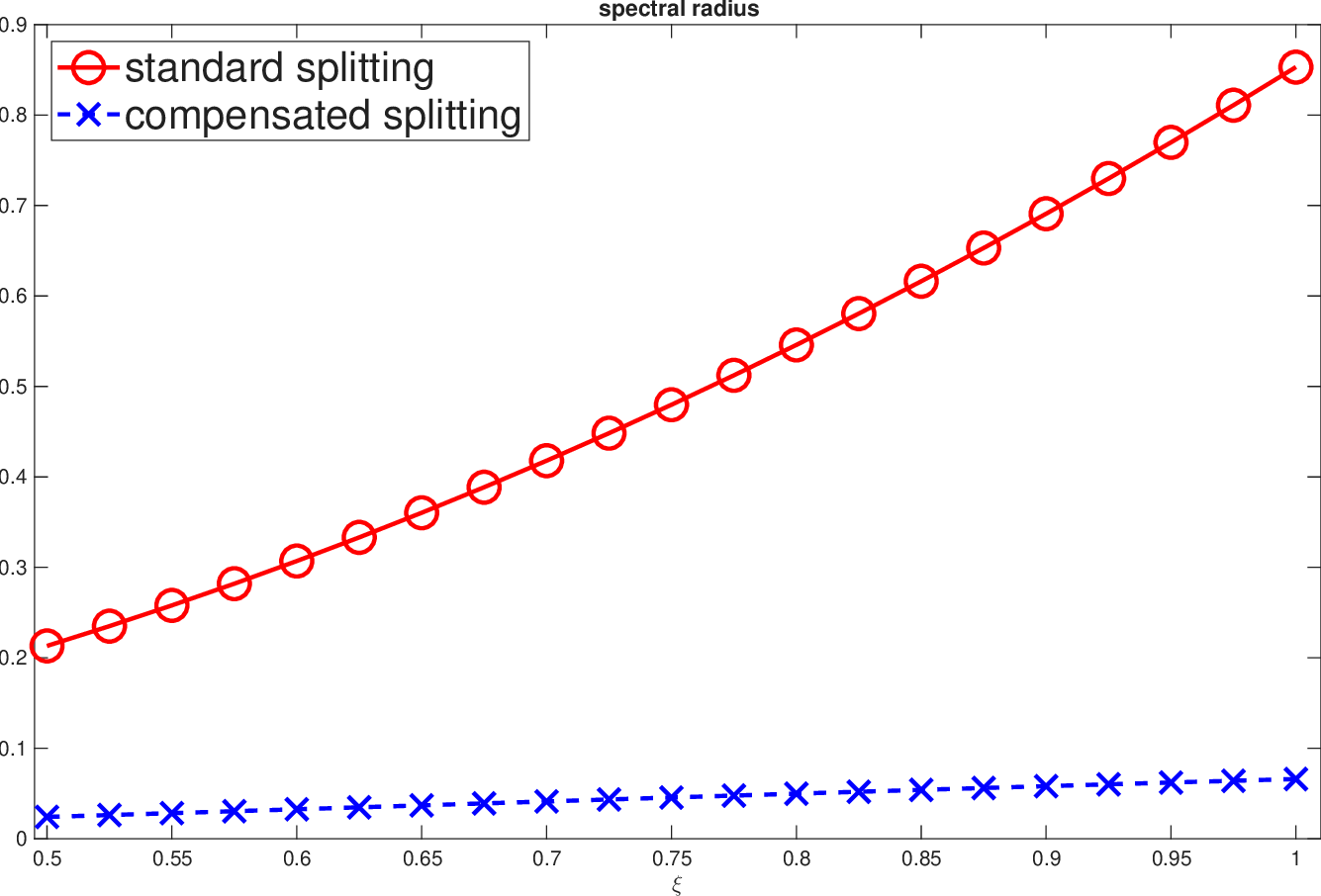}}}\,
\subfloat[number of FP iterations]{\resizebox*{0.32\textwidth}{0.17\textheight}{\includegraphics{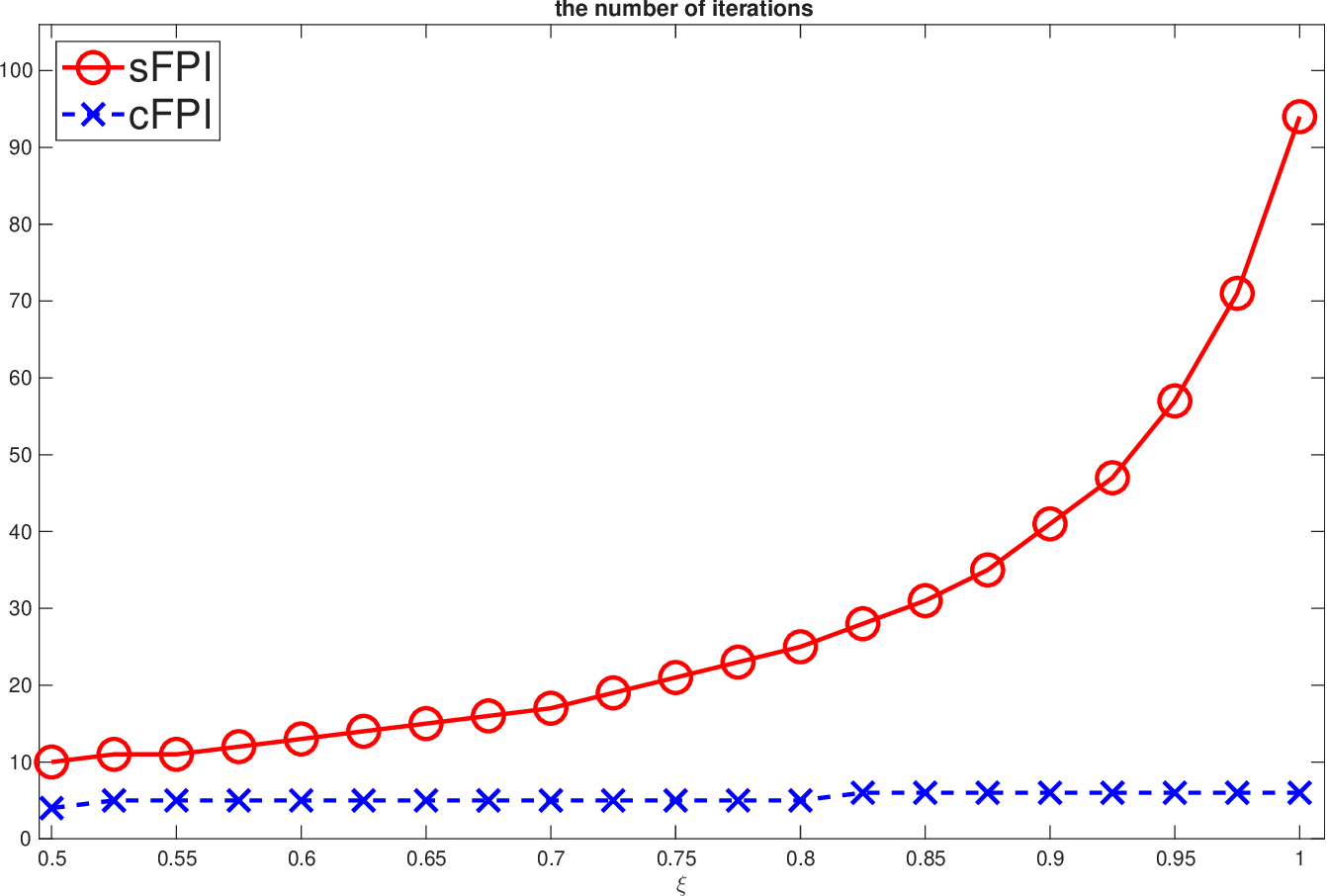}}}\,
\subfloat[CPU in seconds]{\resizebox*{0.32\textwidth}{0.17\textheight}{\includegraphics{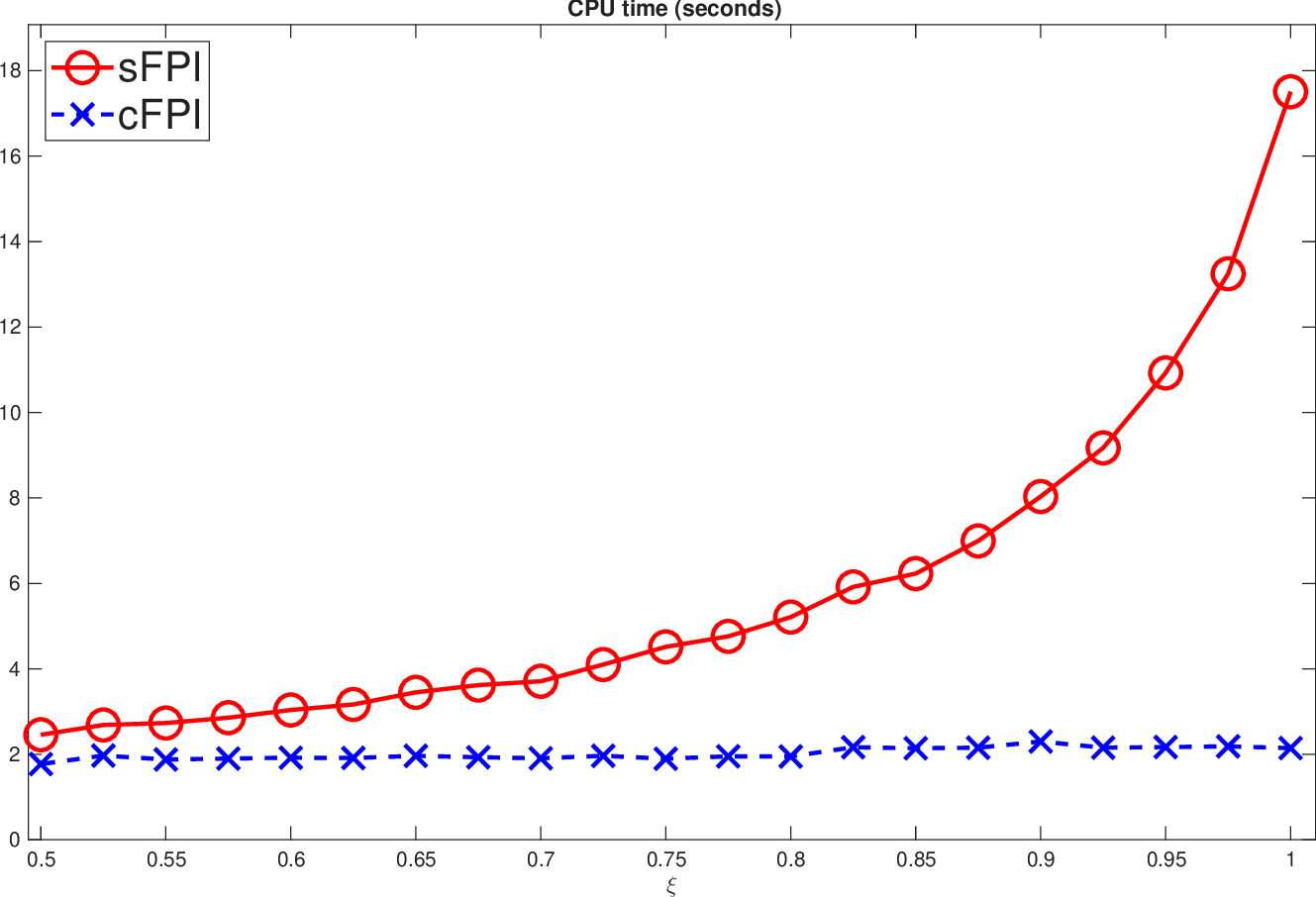}}}
}
\caption{Numerical performance on Example~\ref{eg:eg3} with the first row for $A=A_0$ and the second row for $A=-A_0$,
         as $\xi$ varies in $[0.5,1]$.
  }\label{fig:eg3}
\end{figure}

\begin{example}\label{eg:eg3}
{\rm
The matrix $A$ is the well-known matrix
as in  \texttt{HEAT1} \cite{benner2013low}.
We take $\ell=40$ and define the mesh size $h = 1/(\ell+1)$, and let $\be_1, \be_{\ell}\in\bbR^{\ell}$ be the
first and last columns of $I_{\ell}$,  respectively, and
$$
K =
\left[\begin{array}{ccccc}
2 & -1 & & & \\
-1 & 2 & -1 & & \\
& \ddots & \ddots & \ddots & \\
& & -1 & 2 & -1 \\
& & & -1 & 2
\end{array}\right]\in\bbR^{\ell\times \ell},  \quad
C =
\left[\begin{array}{ccccc}
0 & 1 & & & \\
-1 & 0 & 1 & & \\
& \ddots & \ddots & \ddots & \\
& & -1 & 0 & 1 \\
& & & -1 & 0
\end{array}\right]\in\bbR^{\ell\times \ell},
$$
$E_1 = \be_1\be_1^{\T}$, $E_{\ell} = \be_{\ell} \be_{\ell}^{\T}$. We construct $A_0$ as
\begin{equation}\label{eg3:eq1}
A_0 = \frac{-(I_{\ell}\otimes K + K\otimes I_{\ell})+0.5(E_1\otimes I_{\ell} + E_{\ell}\otimes I_{\ell})}{h^2} + \frac{I_{\ell}\otimes C}{2h}
  \in\bbR^{n\times n}
\end{equation}
with $n=\ell^2$. We take $B = \texttt{randn}(n,5)$, $m=1$ and
\begin{equation*}
N_1=\xi
\begin{bmatrix}
7.8 & -1 & -0.05 &  &  \\
-1 & 7.8 & -1 & \ddots &  \\
-0.05 & -1 & \ddots & \ddots & -0.05 \\
 & \ddots & \ddots & \ddots & -1 \\
 &  & -0.05 & -1 & 7.8
\end{bmatrix},
\end{equation*}
where $\xi$ is a parameter varying from $0.5$ to $1$. Provably,
$\rho(\scrM^{-1}\scrN)=c\xi^2$ where $c>0$ is a constant independent of $\xi$ and thus increases as $\xi$ does.
This fact holds for the two other examples later.

We will solve two parameterized generalized Lyapunov equations \eqref{eq:genlyap} with $A=\pm A_0$, respectively, with
$\ell=40$ yielding $n=1600$.
In \Cref{fig:eg3}, with respect to
the standard splitting \eqref{eq:NatSplit} and our compensated splitting \eqref{eq:t-scrMscrN} with \eqref{eq:optimal(E):intro},
we compare, as $\xi$ varies, spectral radius, the number of iterations, and CPU in seconds, for the case $A=A_0$ (first row)
and $A=-A_0$ (second row). We have the following observations:
\begin{enumerate}[(a)]
  \item For the range of $\xi\in [0.5,1]$, both $\rho(\scrM^{-1}\scrN)$ and $\rho(\wtd\scrM^{-1}\wtd\scrN\,)$ are less than $1$, indicating both
        sFPI and the proposed cFPI are convergent.
  \item It makes no difference in the spectral radius $\rho(\scrM^{-1}\scrN)$ whether $A=A_0$ or $-A_0$. This is expected because
        the spectral radius is invariant with respect to $A=\pm A_0$. But the story is entirely different for the compensated splitting. Besides
        $\rho(\wtd\scrM^{-1}\wtd\scrN\,)$ is smaller than $\rho(\scrM^{-1}\scrN)$ at each $\xi$ in both cases,
        $\rho(\wtd\scrM^{-1}\wtd\scrN\,)$ for the case $A=-A_0$ is much smaller than it for the case $A=A_0$ at the same $\xi$.
        The phenomenon can be heuristically explained by our earlier discussion following \Cref{thm:GCS(complex)}, as
        the eigenvalues of $A_0$ lie in the left-half plane.
  \item For the case $A=-A_0$ whose eigenvalues lie in the right-half plane, $\rho(\wtd\scrM^{-1}\wtd\scrN\,)$ increases very slowly as $\xi$ increases and always stays below $0.1$, while
        $\rho(\scrM^{-1}\scrN)$ grows at a much faster pace and ends up being $0.85$ at $\xi=1$.
        As a result, cFPI runs much faster than sFPI, respectively at each $\xi$.
\end{enumerate}
}
\end{example}

\begin{figure}[t]
{\centering
{\resizebox*{0.32\textwidth}{0.17\textheight}{\includegraphics{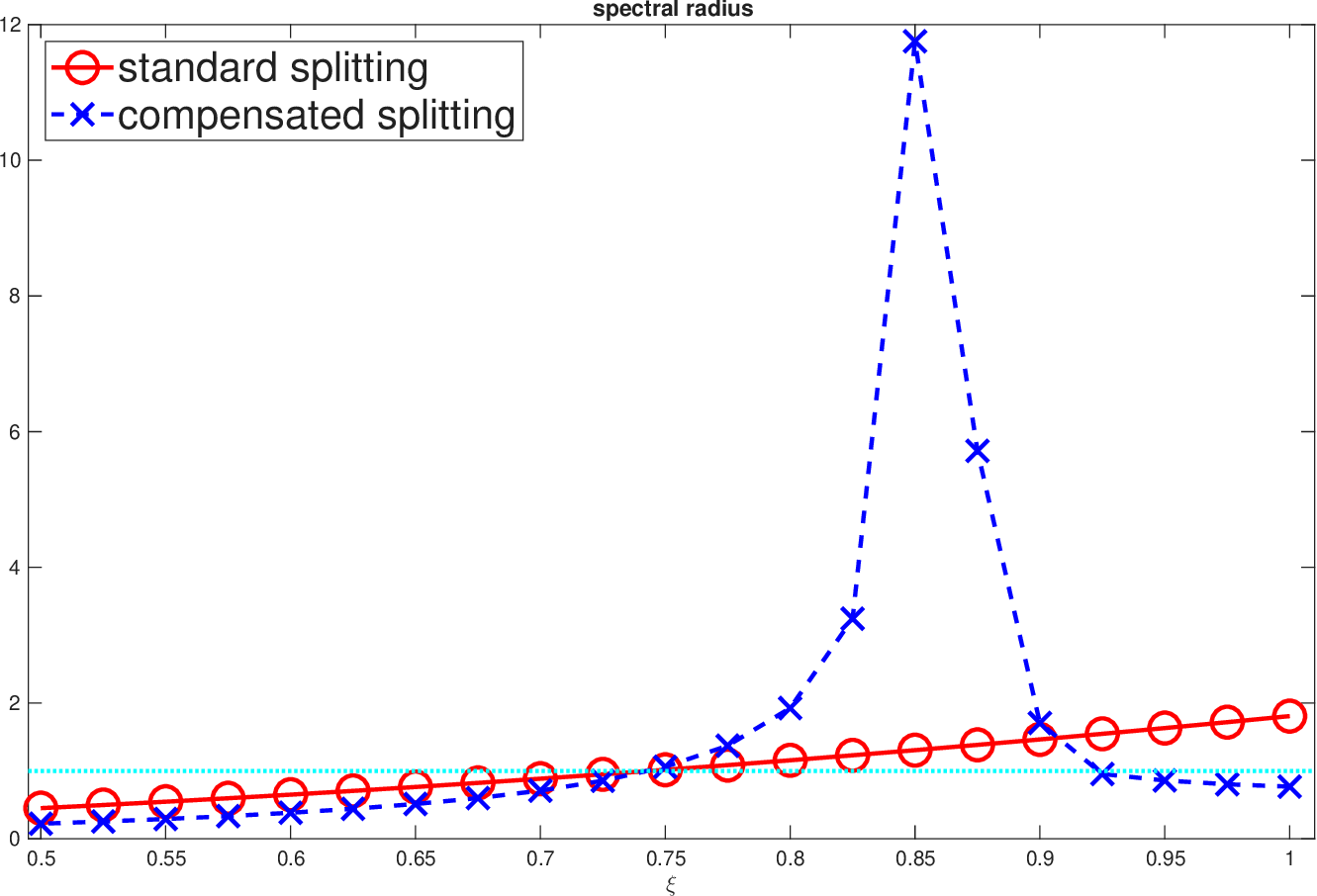}}}\,
{\resizebox*{0.32\textwidth}{0.17\textheight}{\includegraphics{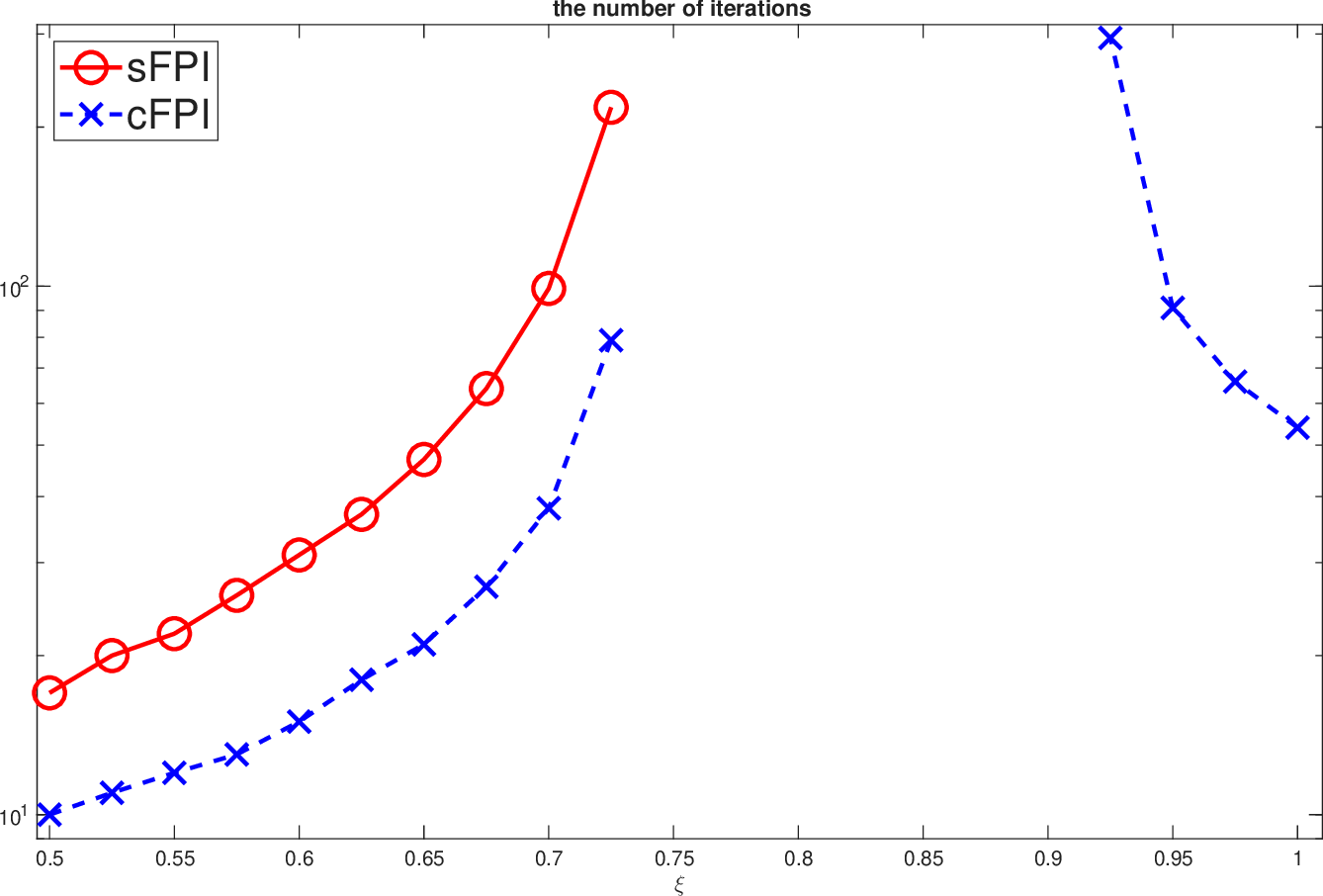}}}\,
{\resizebox*{0.32\textwidth}{0.17\textheight}{\includegraphics{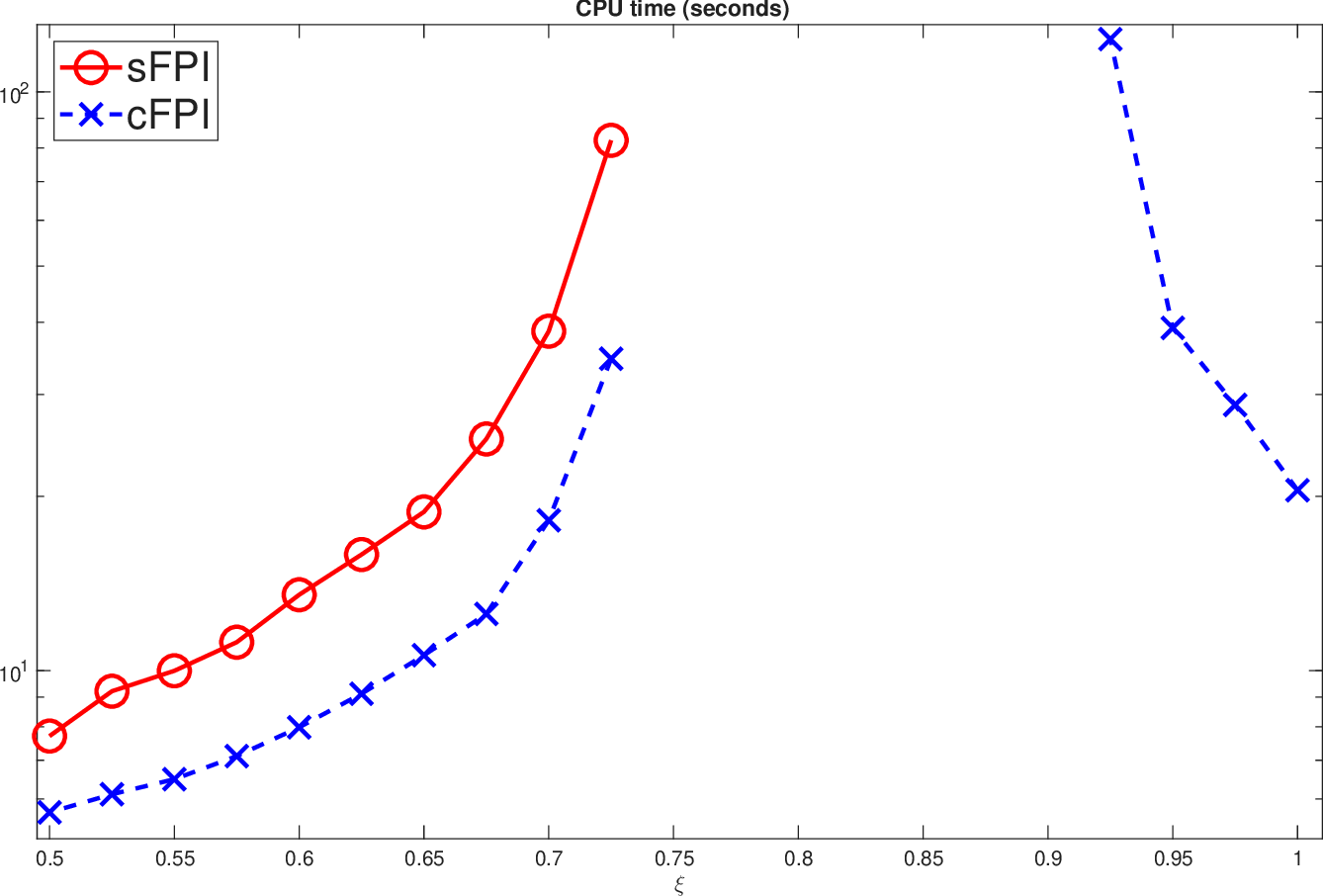}}}
} \\
{\centering
\subfloat[spectral radius]{\resizebox*{0.32\textwidth}{0.17\textheight}{\includegraphics{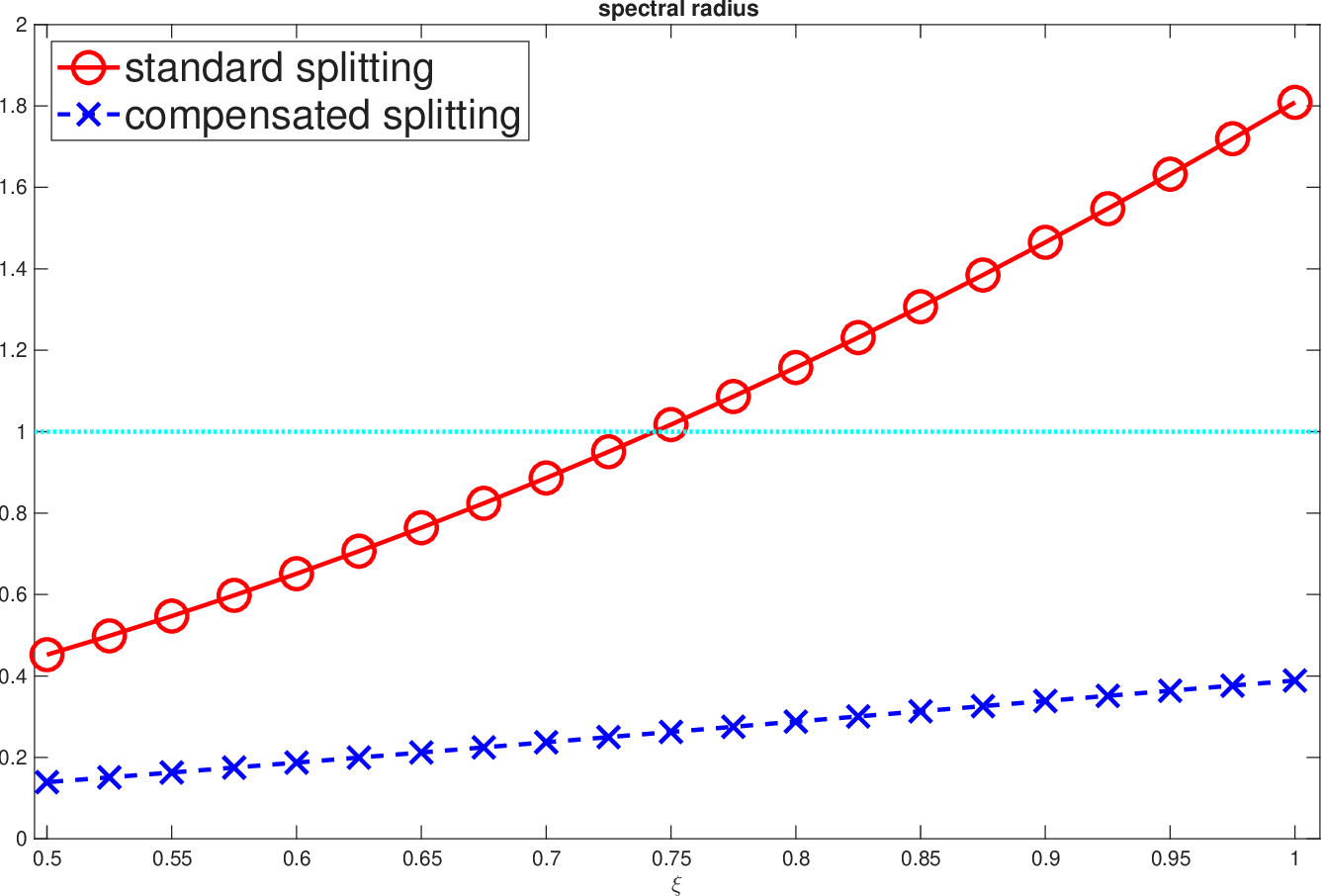}}}\,
\subfloat[number of FP iterations]{\resizebox*{0.32\textwidth}{0.17\textheight}{\includegraphics{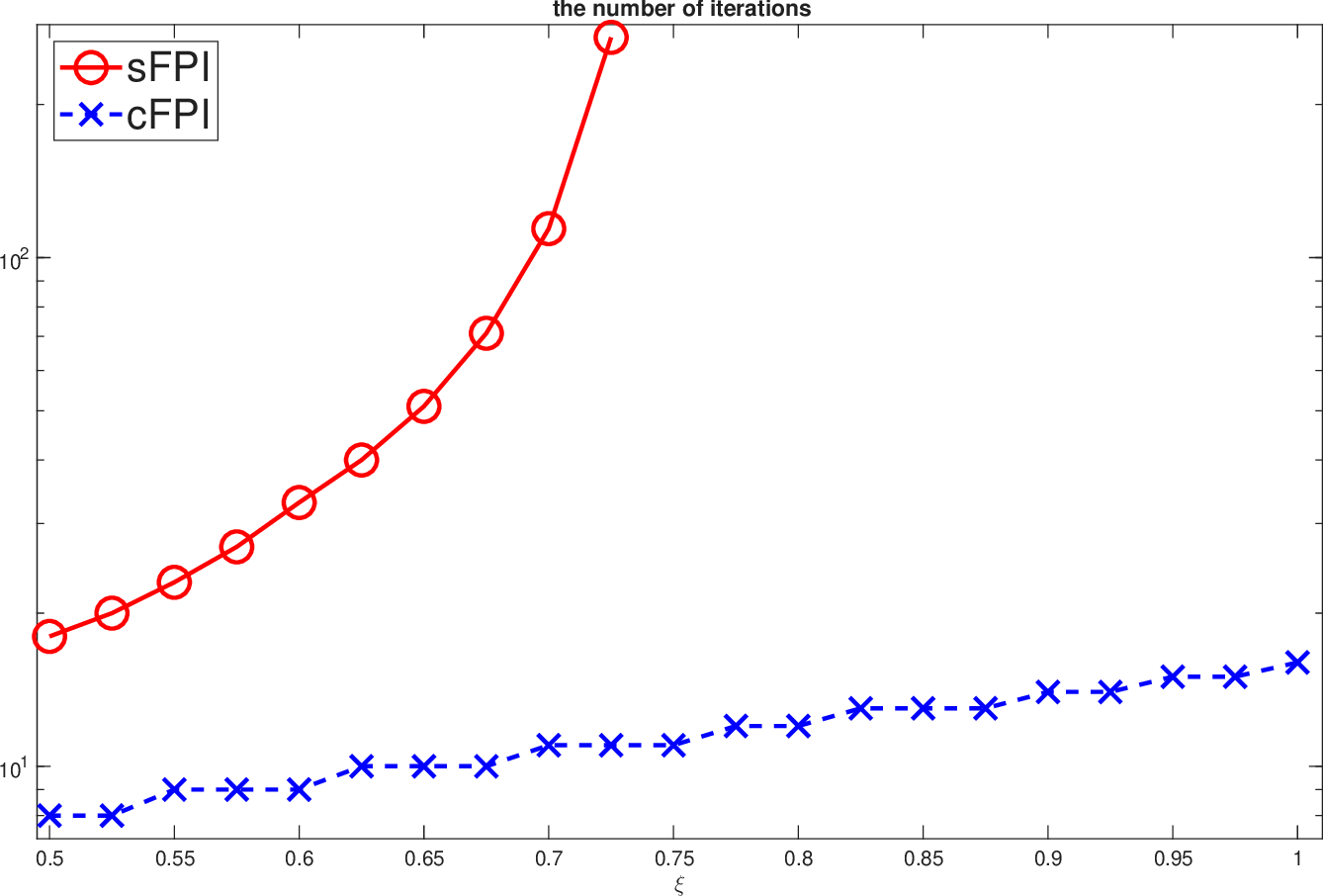}}}\,
\subfloat[CPU in seconds]{\resizebox*{0.32\textwidth}{0.17\textheight}{\includegraphics{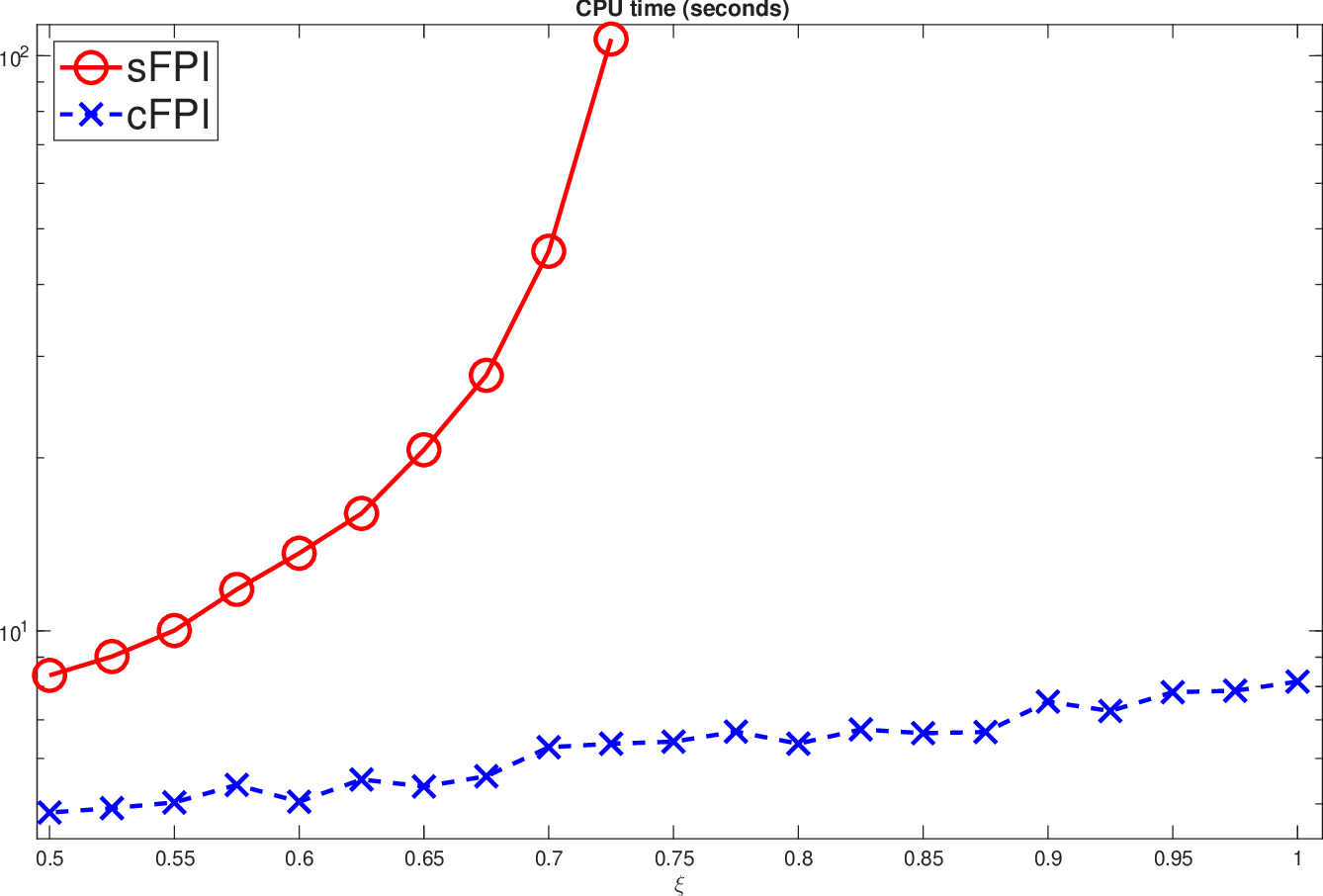}}}
}
\caption{Numerical performance on Example~\ref{eg:egx:r}  with the first row for $A=A_0$ and the second row for $A=-A_0$,
         as $\xi$ varies in $[0.5,1]$.
  No numerical result is shown when the associated fixed-point iteration is deemed divergent according to
  the dotted line for $\rho=1$ in the plots of the left panel.
  }\label{fig:egx:r}
\end{figure}

\begin{example}\label{eg:egx:r}
{\rm
The matrix $A=\pm A_0$with $A_0$ being the tridiagonal matrix
\begin{equation}\label{eq:1-2-1}
A_0 =
\left[\begin{array}{ccccc}
-2 & 1 & & & \\
1 & -2 & 1 & & \\
& \ddots & \ddots & \ddots & \\
& & 1 & -2 & 1 \\
& & & 1 & -2
\end{array}\right]\in\bbR^{n\times n},
\end{equation}
$n=1000$, $B=\texttt{randn}(n,2)\in\bbR^{n \times 2}$,
$m=5$, and $N_j = 2\cdot 10^{-4}\times\xi\times C_jC_j^{\T}$ with each
$C_j=\texttt{randn}(n,10)$, where  $\xi$ is a parameter varying from $.5$ to $1$.
Numerical results are shown in \Cref{fig:egx:r}. We have the following observations:
\begin{enumerate}[(i)]
  \item Provably $\rho(\scrM^{-1}\scrN)$ monotonically increases and it crosses $1$ at $\xi$ somewhere between $0.725$ and $0.75$.
        Because of that,  sFPI converges for $\xi\le 0.725$  and diverges otherwise.
  \item Interestingly, for the case $A=A_0$, $\rho(\wtd\scrM^{-1}\wtd\scrN\,)$ shows a ``pop-up'' starting at $\xi=0.75$ until $\xi=0.9$.
        Excluding the period of the ``pop-up'', it appears $\rho(\wtd\scrM^{-1}\wtd\scrN\,)$ is less than $\rho(\scrM^{-1}\scrN)$.  The sudden ``pop-up''
        may be explained by our earlier discussion following \Cref{thm:GCS(complex)}: $A_0$ is negative definite and its
        eigenvalues lie in the left-half plane, and $E$ collectively moves the eigenvalues of $A=A_0$ rightward towards $0$,
        potentially crossing over to the right-side of $0$.
  \item  Before $\xi=0.75$ and after $\xi=0.9$,  $\rho(\wtd\scrM^{-1}\wtd\scrN\,)$ is smaller than $1$, and thus cFPI converges  for $\xi$ before $0.75$ and after $0.90$.
  \item As in \Cref{eg:eg3}, for the case $A=-A_0$, $\rho(\wtd\scrM^{-1}\wtd\scrN\,)$ increases very slowly as $\xi$ increases.
        It is noted that $A=-A_0$ is positive definite and its eigenvalues lie in the right-half plane.
        Collectively, $E$ moves the eigenvalues of $A=-A_0$ rightward further away from $0$.
\end{enumerate}
}
\end{example}

\begin{figure}[t]
{\centering
{\resizebox*{0.32\textwidth}{0.17\textheight}{\includegraphics{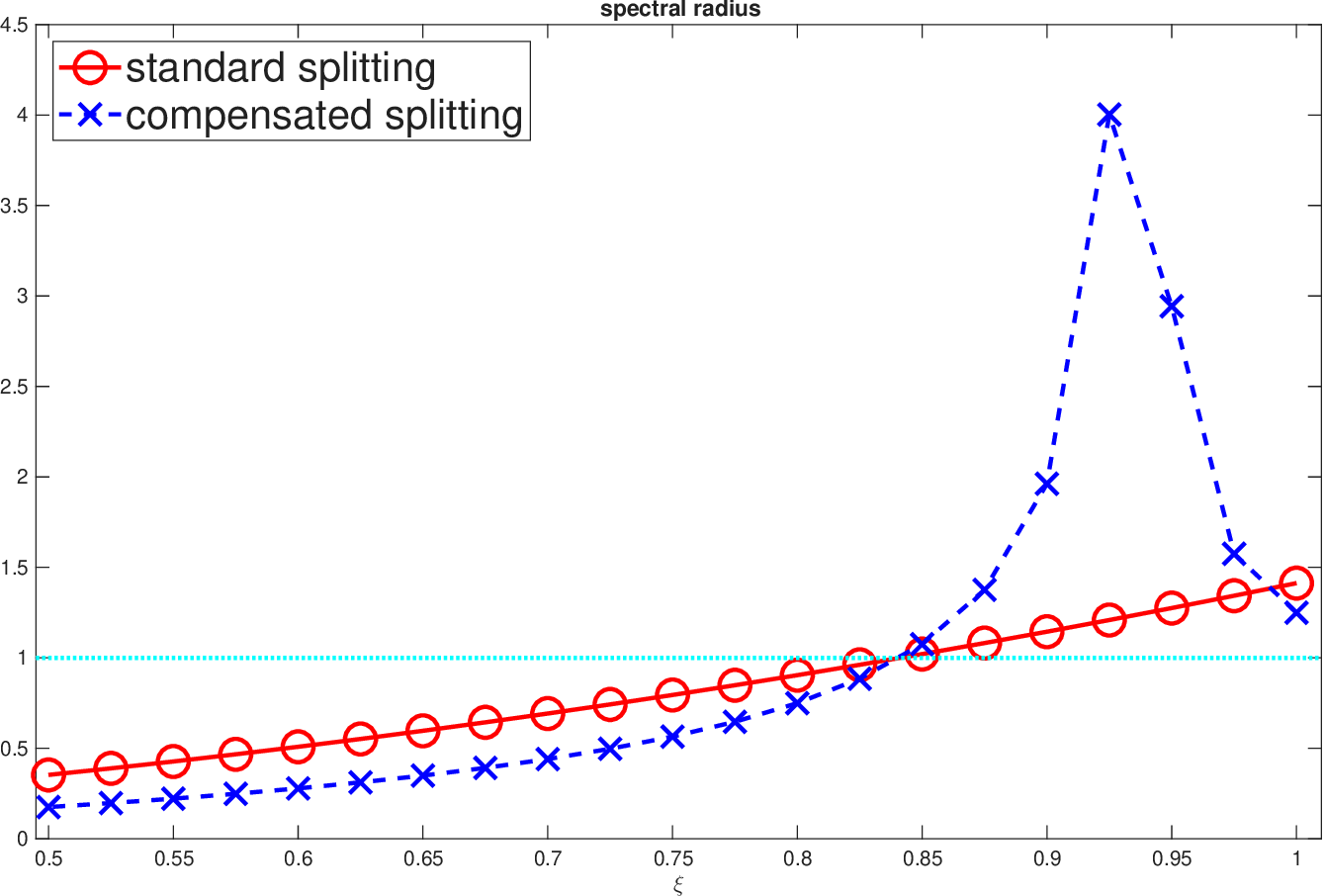}}}\,
{\resizebox*{0.32\textwidth}{0.17\textheight}{\includegraphics{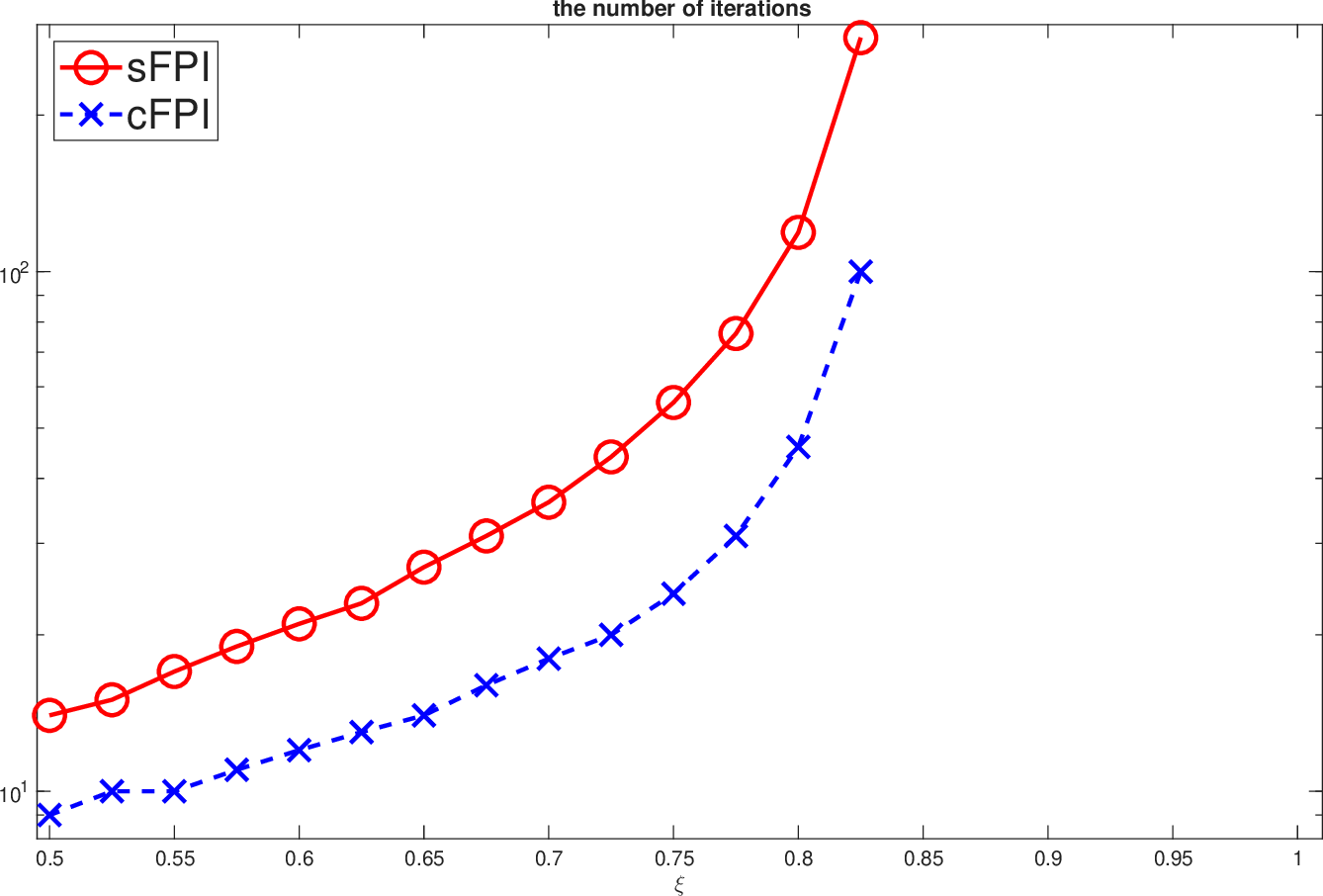}}}\,
{\resizebox*{0.32\textwidth}{0.17\textheight}{\includegraphics{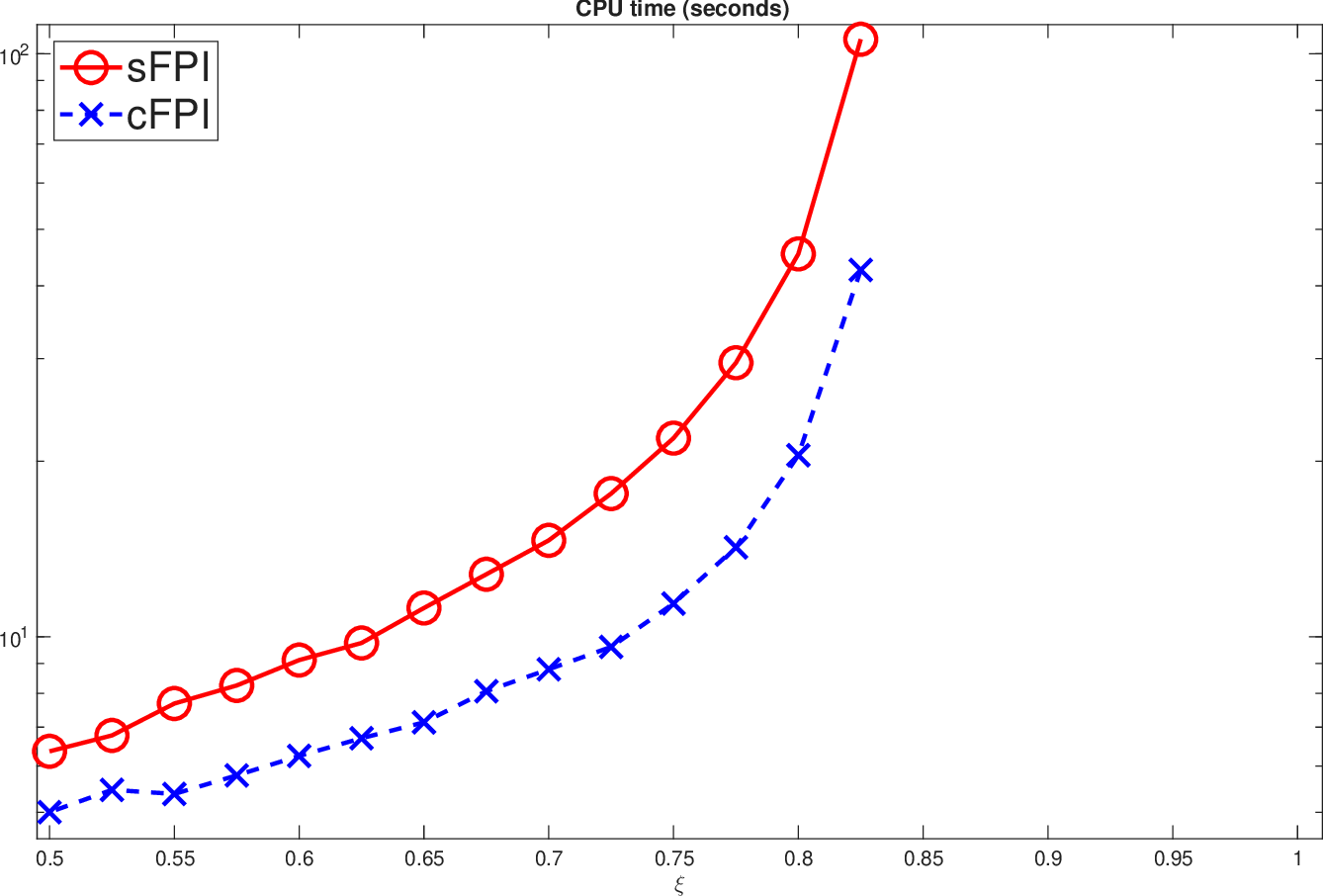}}}
} \\
{\centering
\subfloat[spectral radius]{\resizebox*{0.32\textwidth}{0.17\textheight}{\includegraphics{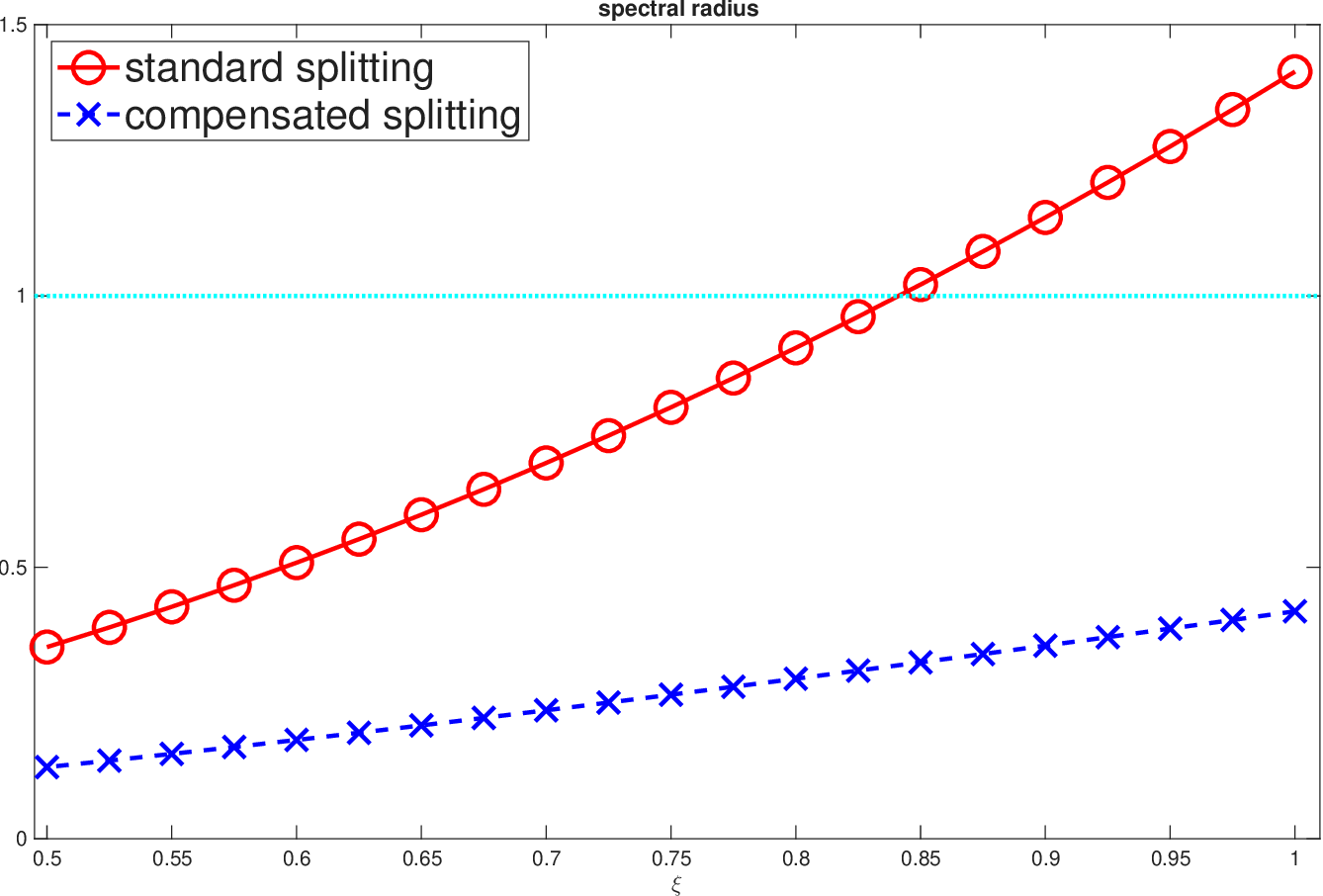}}}\,
\subfloat[number of FP iterations]{\resizebox*{0.32\textwidth}{0.17\textheight}{\includegraphics{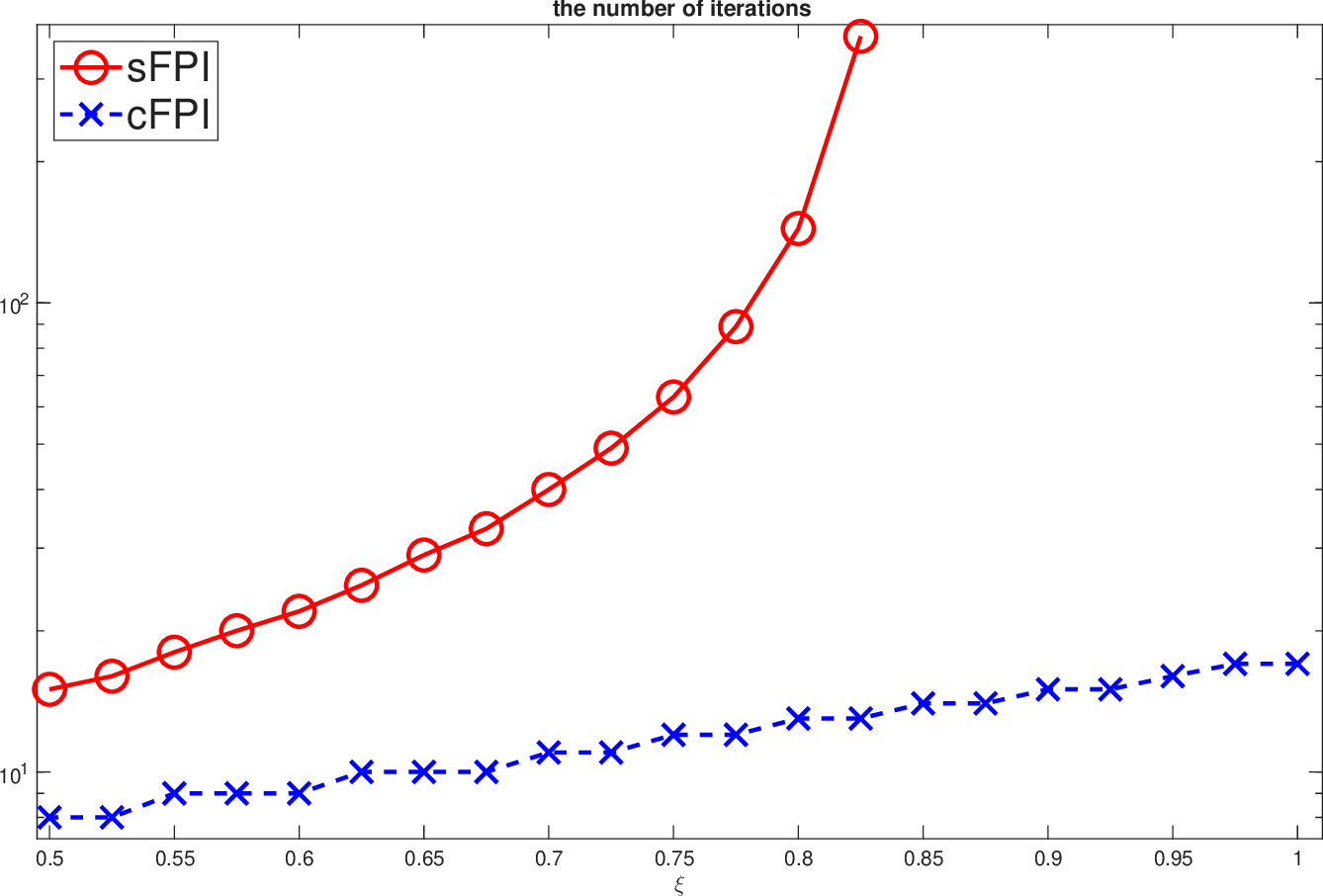}}}\,
\subfloat[CPU in seconds]{\resizebox*{0.32\textwidth}{0.17\textheight}{\includegraphics{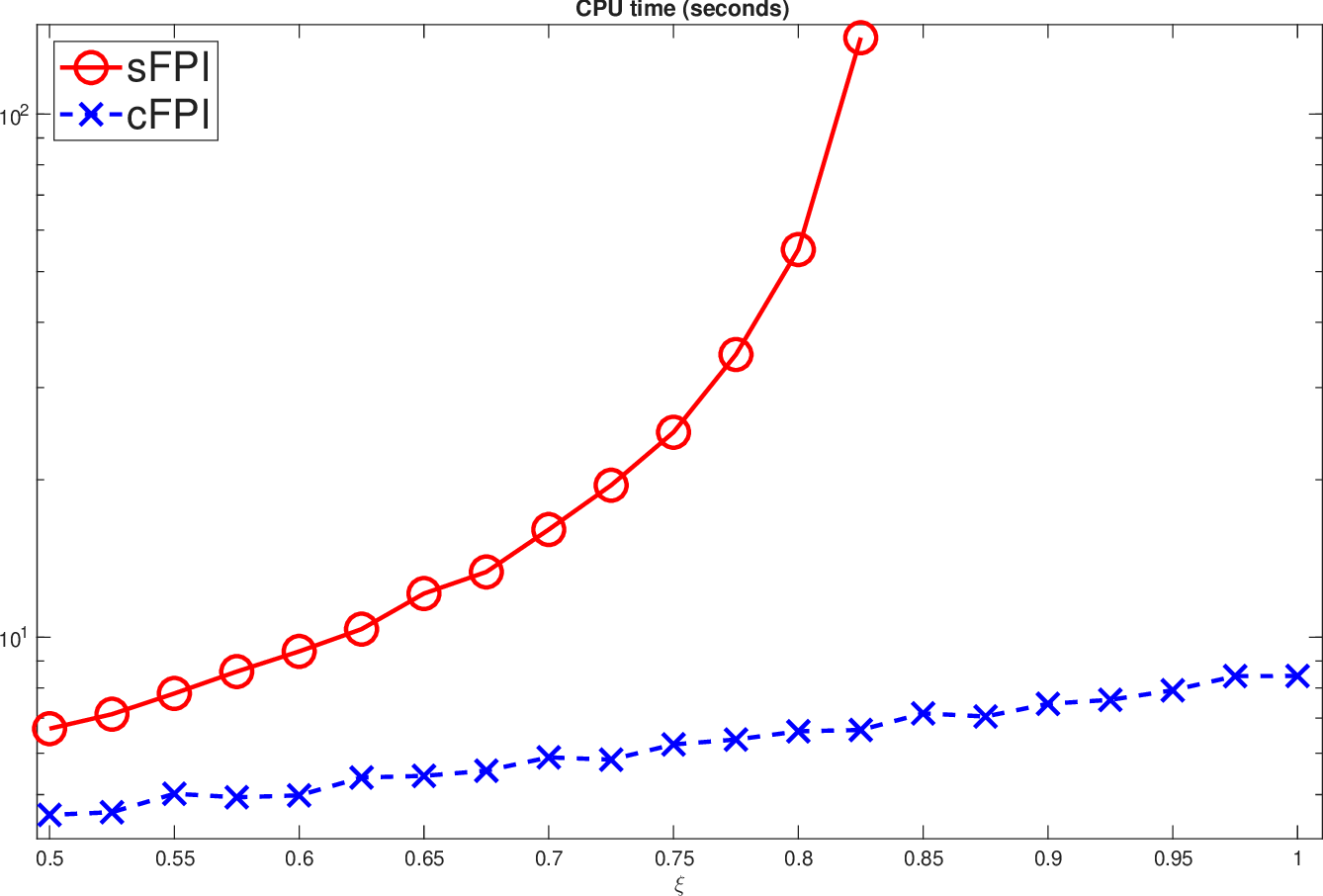}}}
}
\caption{Numerical performance on Example~\ref{eg:egx:c}  with the first row for $A=A_0$ and the second row for $A=-A_0$,
         as $\xi$ varies in $[0.5,1]$.
  No numerical result is shown when the associated fixed-point iteration is deemed divergent according to
  the dotted line for $\rho=1$ in the plots of the left panel.
  }\label{fig:egx:c}
\end{figure}

\begin{example}\label{eg:egx:c}
{\rm
This is a complex example. We add a small imaginary number to the diagonal of the one in \eqref{eq:1-2-1} to get
\begin{equation}\label{eq:1-2-1:c}
A_0=
\left[\begin{array}{ccccc}
-2 & 1 & & & \\
1 & -2 & 1 & & \\
& \ddots & \ddots & \ddots & \\
& & 1 & -2 & 1 \\
& & & 1 & -2
\end{array}\right]+{\tt 1i}\times 10^{-8}\times I_n\in\bbC^{n\times n},
\end{equation}
$n=1000$,
$B=\texttt{randn}(n,2)+\texttt{1i*randn}(n,2)$,
$m=5$, and $N_j = 2\cdot 10^{-4}\times\xi\times C_jC_j^{\HH}$ with each
$$
C_j=\texttt{randn}(n,5)+\texttt{1i*randn}(n,5),
$$
where  $\xi$ is a parameter varying from $.5$ to $1$.
Adding a small imaginary number in \eqref{eq:1-2-1:c} is to ensure that calling the Schur decomposition of $A$
at the beginning of sFPI yields a complex Schur form. This will equalize the computational complexity
per fixed-point iterative steps in both sFPI and cFPI.
Numerical results are shown in \Cref{fig:egx:c}. We witness similar things to \Cref{fig:egx:r}:
\begin{enumerate}[(1)]
  \item Provably $\rho(\scrM^{-1}\scrN)$ monotonically increases and it crosses $1$ at $\xi$ somewhere between $0.825$ and $0.85$.
        Because of that,  sFPI converges for $\xi\le 0.825$   and diverges otherwise.
  \item For the case $A=A_0$, $\rho(\wtd\scrM^{-1}\wtd\scrN\,)$ shows a ``pop-up'' starting at $\xi=0.85$ until likely some $\xi>1$ if we let it go larger. Once again, the sudden ``pop-up''
        may be explained by our earlier discussion following \Cref{thm:GCS(complex)}, as the
        eigenvalues of $A=A_0$ lie in the left-half plane.
  \item As in the previous two examples, for the case $A=-A_0$, $\rho(\wtd\scrM^{-1}\wtd\scrN\,)$ increases very slowly as $\xi$ increases.
        It is verified numerically that $A=-A_0$ has its eigenvalues lying in the right-half plane.
\end{enumerate}
}
\end{example}

\section{Concluding Remarks}\label{sec:conclude}
The generalized Lyapunov equation
$$
AX + XA^{\HH} + \sum_{j=1}^{m}N_j X N_j^{\HH} + BB^{\HH} = 0
$$
arises frequently in applications.
Such an equation is naturally split into $\scrM(X)=\scrN(X)-BB^{\HH}$ where $\scrM(X)=AX + XA^{\HH}$ and
$\scrN(X)=-\sum_{j=1}^{m}N_j X N_j^{\HH}$, and it has been done frequently to yield the commonly used
standard fixed-point iteration (sFPI) for solving the equation.
The iteration converges
if the spectral radius $\rho(\scrM^{-1}\scrN)<1$ but otherwise diverges in general.
In this paper, we propose
a compensated splitting technique to tackle situations where $\rho(\scrM^{-1}\scrN)\ge 1$
or $\rho(\scrM^{-1}\scrN)<1$ but too close to $1$ leading very slow convergence.
The basic idea is to
compensate $\scrM(X)$ by projecting $\sum_{j=1}^{m}N_j X N_j^{\HH}$ into an operator in the form $EX + XE^{\HH}$ to merge with
$\scrM(X)$, yielding the so-called compensated splitting $\wtd\scrM(X)=\wtd\scrN(X)-BB^{\HH}$ for which there is a
corresponding fixed-point iteration called the compensated fixed-point iteration (cFPI). Ideally,
we should minimize $\rho(\wtd\scrM^{-1}\wtd\scrN\,)$ over matrix $E$, but that is not practical. Instead
we determine an optimal $E$ by \eqref{eq:SARE-CS''}.
It is shown that \eqref{eq:SARE-CS''} has a unique global minimizer in the real case for every $\mu\ge0$.
In the complex case, the global minimizer is unique when $\mu>0$;
when $\mu=0$, there are infinitely many global minimizerglobal minimizers but we select the one with the smallest Frobenius norm.
Compactly, the optimal $E$ is given by
$$
E=\frac{2}{2n+\mu}\sum_{j=1}^m\tr(\bN_{j})\left[N_j-\frac{2}{4n+\mu}\tr(N_{j})\cdot I_n\right],
$$
where $\mu\ge 0$ is a regularizing parameter to control the magnitude of $E$. At this point, we have little idea on
how to determine a good $\mu$, and in fact in our reported numerical experiments in \cref{sec:egs}, we simply take
$\mu=0$ as the default.

Unfortunately, our current \eqref{eq:SARE-CS''} to determine $E$ does not involve the $\scrM$-part in any way. Heuristically,
we have analyse the interaction between this optimal $E$ and the $\scrM$-part. Our numerical experiments
appear to behave consistently with our analysis, while
demonstrating  superiority of cFPI to sFPI.

In our numerical experiments, we limit the problem sizes, i.e., $n$, modest so that
involved Lyapunov equations in the middle of computations can be solved by the Bartels-Stewart method \cite{bast:1972}
with ease on everyday personal computers. For that reason,
our purpose in this paper is two-fold: one is to demonstrate our compensated splitting can
rescue the natural splitting, the commonly used one \eqref{eq:NatSplit} when the resulting
sFPI \eqref{eq:NatItn} diverges, and the other one is to show that the compensated splitting
can speed up cFPI when it indeed converges. It is expected the key idea in this paper can be combined with
more sophisticated subspace-based projection iterative methods, e.g., \cite{breiten2019residual,shank2016efficient},
to yield more efficient methods for large scale generalized Lyapunov equations.
Such investigations will be forthcoming.

A similar looking to the generalized Lyapunov equation is the generalized Sylvester equation
\cite{jamp:2018,krsi:2015,simo:2016,voet:2025}
$$
AX + XB^{\HH} + \sum_{j=1}^{m}N_j X M_j^{\HH} + C = 0.
$$
It arises from, for example,  discretizing certain partial differential equations \cite{pasi:2016} and perhaps other places, and it too
possesses a natural splitting that leads to a natural fixed-point iteration, similarly to the discussions
we had in \cref{sec:intro} for the generalized Lyapunov equation. Due to high similarity, our compensated splitting
in \cref{sec:splitting} can be straightforwardly attempted for potential improvement. We will leave the detail to another
paper.


{\small
\def\noopsort#1{}\def\l{\char32l}\def\v#1{{\accent20 #1}}
  \let\^^_=\v\def\hbk{hardback}\def\pbk{paperback}

}

\end{document}